\documentclass[11pt]{article}
\usepackage{amssymb,epsfig,amsmath}
\usepackage[dvips]{color}
\usepackage[T1]{fontenc}
\usepackage[francais]{babel}
\usepackage[utf8]{inputenc}
\usepackage{comment}
\usepackage{hyperref}
\usepackage{graphicx}
\usepackage{multicol}
\usepackage{caption}
\newtheorem{defi}{D\'efinition}[section]
\newtheorem{prop}[defi]{Proposition}
\newtheorem{theo}[defi]{Th\'eor\`eme}
\newtheorem{conj}[defi]{Conjecture}
\newtheorem{lemm}[defi]{Lemme}
\newtheorem{coro}[defi]{Corollaire}
\newtheorem{rema}[defi]{Remarque}
\newtheorem{exem}[defi]{Exemple}
\newtheorem{exems}[defi]{Exemples}

\newcommand{\bdefi}{\begin{defi}}
\newcommand{\edefi}{\end{defi}}
\newcommand{\bprop}{\begin{prop}}
\newcommand{\eprop}{\end{prop}}
\newcommand{\btheo}{\begin{theo}}
\newcommand{\etheo}{\end{theo}}
\newcommand{\blemm}{\begin{lemm}}
\newcommand{\brema}{\begin{rema}}
\newcommand{\erema}{\end{rema}}
\newcommand{\bexer}{\begin{exem}}
\newcommand{\eexer}{\end{exem}}
\newcommand{\bexems}{\begin{exems}}
\newcommand{\eexems}{\end{exems}}
\newcommand{\bconj}{\begin{conj}}
\newcommand{\econj}{\end{conj}}
\newcommand{\elemm}{\end{lemm}}
\newcommand{\bcoro}{\begin{coro}}
\newcommand{\ecoro}{\end{coro}}
\newcommand{\dem}{\noindent{\bf Démonstration. }}
\newcommand{\rem}{\noindent{\bf Remarque. }}

\usepackage{mathrsfs}
\renewcommand\mathcal{\mathscr}

\newcommand{\G}{{\cal G}}

\newcommand{\F}{{\cal F}}

\newcommand{\C}{{\cal C}}

\newcommand{\Q}{{\cal Q}}

\newcommand{\maths}[1]{{\mathbb #1}}  

\newcommand{\EE}{\maths{E}}

\newcommand{\RR}{\maths{R}}
\newcommand{\NN}{\maths{N}}
\newcommand{\CC}{\maths{C}}

\newcommand{\SSS}{\maths{S}}

\newcommand{\Ga}{\Gamma}

\newcommand{\cqfd}{\hfill$\Box$}

\newcommand{\card}{{\operatorname{Card}}}

\renewcommand{\Im}{{\operatorname{Im}}}

\newcommand{\CAT}{\operatorname{CAT}}

\newcommand{\dddp}{\partial_{\infty}^2\widetilde{\Sigma}}
\newcommand{\ddp}{\partial_{\infty}\widetilde{\Sigma}}
\newcommand{\ddX}{\partial_{\infty}\widetilde{X}}

\newcommand{\ddT}{\partial_\infty T}

\newcommand{\Plat}{\operatorname{Plat}}

\newcommand{\Diff}{\operatorname{Diff}}

\newcommand{\revet}{(\widetilde{\Sigma},[\widetilde{q}])}

\newcommand{\srfce}{({\Sigma},{[q]})}

\newcommand{\revetm}{(\widetilde{\Sigma},\widetilde{m})}
\newcommand{\srfcem}{({\Sigma},m)}
\newcommand{\Id}{\operatorname{Id}}
\newcommand{\pl}{\operatorname{pl}}
\newcommand{\Lq}{\Lambda_{[q]}}

\newcommand{\Lqr}{\widetilde{\Lambda}_{[\widetilde{q}]}}
\newcommand{\Lm}{\Lambda_m}
\newcommand{\Lmr}{\widetilde{\Lambda}_{\widetilde{m}}}
\title{Laminations géodésiques plates.}
\author{T.~Morzadec} 

\begin{document}
\maketitle

\textbf{Abstract:} Since their introduction by Thurston, geodesic laminations on hyperbolic surfaces occur in many contexts.
In this paper, we propose a generalization of geodesic laminations on locally $\CAT(0)$, complete, geodesic metric spaces, whose boundary
at infinity of the universal cover
is endowed with an  invariant total cyclic order. Then we study these new objects on surfaces endowed with half-translation structures 
and on finite metric graphs. The main result of the paper is a theorem of classification of geodesic laminations on a compact surface
endowed with a half-translation structure. We also show that every finite metric fat graph, outside four homeomorphism classes, is the support of a geodesic lamination with uncountably many
leaves none of which is eventually periodic. \footnote{ Mots clés: Lamination géodésique, surface munie d'une structure de demi-translation, 
    différentielle quadratique holomorphe, 
feuilletage singulier, surface hyperbolique.
Code AMS 30,37,53,57,58.}

\section{Introduction.}
L'objet principal de cet article est de proposer et d'étudier une généralisation des laminations géodésiques sur les surfaces hyperboliques (voir par exemple 
\cite{Bonahon97}) à celles sur les espaces métriques localement $\CAT(0)$ enrubannés (voir ci-dessous), plus particulièrement aux surfaces munies de
structures de demi-translation
(c'est-à-dire de structures plates à 
holonomie $\{\pm\Id\}$). Une motivation, qui fera l'objet d'un travail ultérieur, est d'étudier en profondeur les dégénérescences de tels objets, comme initié
dans \cite{DucLeiRaf10}. Dans le cas des arbres, cette notion de lamination géodésique ne doit pas être confondue avec celle introduite dans par exemple
\cite{BesFeiHan97} et \cite{CouHilLus07}.

Soient $(X,d)$ un espace métrique géodésique complet, localement $\CAT(0)$, et  $p:(\widetilde{X},\widetilde{d})\to(X,d)$
  un rev\^etement universel localement isométrique. On note $\partial_\infty\widetilde{X}$ 
le bord à l'infini de $(\widetilde{X},\widetilde{d})$ et $\Gamma$ le groupe de rev\^etement de $p$. 
Un \textit{ordre cyclique} sur $\partial_\infty\widetilde{X}$ est une application
$o:(\partial_\infty\widetilde{X})^3\to\{-1,0,1\}$, invariante pour l'action 
diagonale de $\Gamma$ sur $(\partial_\infty\widetilde{X})^3$, analogue à l'ordre cyclique sur le cercle défini par un choix d'orientation 
(c'est-à-dire $o(x,y,z)$ vaut $0$
si $x,y$ et $z$ ne sont pas deux à deux distincts, vaut $1$ si $x,y$ et $z$ sont ordonnés dans le sens direct, et $-1$ sinon), voir par exemple \cite[2.3.1]{Wolf11}.

%
%
%

On dit que $(X,d)$ est \textit{enrubanné} si $\partial_\infty\widetilde{X}$ est muni d'un tel ordre cyclique $o$.
Deux couples de points $(x_1,y_1)$ et 
$(x_2,y_2)$ de $\partial_\infty\widetilde{X}$ sont dits \textit{entrelacés} si les points sont deux à deux distincts et si $o(x_1,x_2,y_1)=-o(x_1,y_2,y_1)$, 
et deux géodésiques locales de $(X,d)$ sont dites entrelacées si elles ont des relevés dans
$\widetilde{X}$ dont les couples de points à l'infini sont entrelacés.  
On munit l'ensemble des géodésiques locales, orientées mais non paramétrées, de $(X,d)$ de la topologie quotient de la topologie compacte-ouverte par 
l'action par translations à la source de $\RR$ sur les géodésiques locales paramétrées, que l'on appelle \textit{topologie des géodésiques}. 
 
\bdefi
Une lamination géodésique de $(X,d,o)$ est 
un ensemble non vide $\Lambda$ de géodésiques locales complètes de $(X,d)$, définies à changements d'origine près, dont les éléments sont appellés feuilles, tel que: 

\medskip

\noindent
$\bullet$~ les feuilles de $\Lambda$ sont non auto-entrelacées et deux à deux non entrelacées;

\noindent
$\bullet$~ $\Lambda$ est invariant par changements d'orientation des feuilles;

\noindent
$\bullet$~ $\Lambda$ est fermé pour la topologie des géodésiques.

\noindent
On appellera \textit{support} de $\Lambda$ la réunion des images des feuilles de $\Lambda$.
\edefi

Soit $\Sigma$ une surface connexe, orientable (sans bord pour simplifier dans cette introduction). Une {\it structure de demi-translation} (ou {\it structure plate à holonomie 
$\{\pm\Id\}$}) 
sur  $\Sigma$ est
la donn\'ee d'une partie discrète $Z$ de $\Sigma$
(\'eventuellement vide) et d'une m\'etrique euclidienne sur $\Sigma-Z$ \`a
singularit\'es coniques d'angles $k\pi$, avec $k\in\NN$ et $k\geqslant 3$ aux points de
$Z$, telle que l'holonomie de tout lacet fermé, $\C^1$ par morceaux,  de $\Sigma-Z$ soit contenue dans $\{\pm\Id\}$. Nous renvoyons 
aux parties \ref{structuresplates} et \ref{diffquadholo} notamment pour le cas à bord.
La surface $\Sigma$ munie d'une structure de demi-translation est alors un espace métrique $(\Sigma,d)$, complet et localement 
$\CAT(0)$, qui est naturellement enrubanné (voir la remarque \ref{ordresurface}). Le résultat principal de cet article est le théorème suivant
de classification des laminations géodésiques sur $(\Sigma,d)$, si $\Sigma$ est compacte. 
 
\btheo
Soit $\Lambda$ une lamination géodésique sur une surface compacte connexe munie d'une structure de demi-translation. Alors $\Lambda$
est la réunion d'un nombre fini de composantes cylindriques, de
composantes minimales (de type récurrent ou graphe fini ou paire de feuilles périodiques opposées)
et de feuilles isolées (pour la topologie des géodésiques) dont
chacun des
bouts aboutit dans une composante minimale ou dans une composante cylindrique.
\etheo

Dans ce théorème, une composante cylindrique est un ensemble maximal de feuilles de $\Lambda$ dont les images sont contenues dans un cylindre 
plat non dégénéré (elles sont alors périodiques), une composante minimale est une sous-lamination  
qui est l'adhérence, pour la topologie des géodésiques, de la réunion d'une feuille $\ell$ et de son inverse. Si $\ell$ est régulière 
(c'est-à-dire ne rencontre pas de singularité), alors la composante est de type récurrent (toutes ses feuilles sont alors d'image dense 
dans un domaine de $\Sigma$, c'est-à-dire l'adhérence d'un ouvert connexe, bordé par des géodésiques locales périodiques, s'il n'est pas égal à $\Sigma$). 
Si l'image de $\ell$ est un graphe fini, et si $\ell$ ou son inverse n'est pas périodique
à partir d'un certain temps, alors la composante est de type graphe fini (toutes
ses feuilles ont alors la même image, et aucune n'est périodique à partir d'un certain temps). On dit qu'un bout d'une géodésique \textit{aboutit}
dans une composante minimale ou dans une composante cylindrique si elle admet un rayon géodésique dont l'image est contenue dans le support de la composante minimale ou
dans le bord du cylindre plat correspondant. Nous renvoyons à la partie \ref{laminationplate} pour plus de détails. 

\medskip

Dans la première partie, nous rappellons les propriétés principales des espaces métriques localement $\CAT(0)$ enrubannés et des 
surfaces munies d'une structure de demi-translation. Dans la deuxième, nous montrons que si un rayon d'une géodésique 
locale non auto-entrelacée d'une surface compacte munie d'une structure de demi-translation n'est pas d'image compacte, il admet 
un sous-rayon d'image dense dans un domaine de $\Sigma$. Puis dans la troisième,
nous établissons une correspondance naturelle (mais non bijective, ce qui nécessite un travail de fond) entre les laminations pour une structure de demi-translation 
et celles pour une
métrique 
hyperbolique complète. Dans la quatrième, nous montrons que tous les graphes enrubannés, connexes, finis, sans 
sommet terminal, sauf quatre à homéomorphismes près, sont le support de laminations minimales, sur une infinité de surfaces munies de structures de demi-translation.
Dans la dernière partie, nous montrons le théorème de structure des laminations géodésiques sur les surfaces compactes munies de structures de demi-translation. 
Plusieurs phénomènes radicalement nouveaux apparaissent par rapport aux laminations hyperboliques: les feuilles d'une lamination plate ne sont généralement pas 
deux à deux disjointes, les laminations plates ne sont pas déterminées par 
leur support (des familles non dénombrables ont le même support), les composantes cylindriques peuvent contenir des familles non dénombrables 
de feuilles. Enfin, il y a trois types de composantes minimales d'une lamination plate sur une surface compacte (feuille périodique parcourue dans les deux sens,
composante minimale de type récurrent ou de type graphe fini).   

\medskip

Je tiens à remercier chaleureusement Frédéric Paulin pour sa disponibilité et ses relectures attentives des versions successives de cet article, ainsi que Yves Benoist,
Sylvain Crovisier,
Bertrand Deroin,
Samuel Lelièvre,
Sara Maloni et
 Pierre Pansu  pour leurs réponses à mes nombreuses questions. 
Enfin, je remercie le rapporteur pour ses très 
nombreuses remarques qui m'ont permis d'améliorer considérablement la rédaction de cet article.

%
%
%
%
%
%
%
%
%
%
%
%
%
%
%
%
%
%
%
%
%
%
%
%
%
%
%
%
%
%
%
%

%
%

\section{Définitions.}

\subsection{Espaces localement $\CAT(0)$.}\label{CAT}

Dans tout cet article, nous utiliserons la terminologie de \cite{BriHae99} concernant les espaces métriques : espaces $\CAT(0)$, $\delta$-hyperbolique... 
En particulier, une {\it géodésique} (resp. une {\it géodésique locale}) 
  d'un espace métrique $(X,d)$ est une application isométrique (resp. localement isométrique) $\ell:I\to X$ où $I$ est un intervalle de $\RR$. Nous dirons 
  {\it segment}, {\it rayon} ou {\it droite géodésique} de $X$ si $I$ est respectivement un intervalle compact, une demi-droite fermée (généralement $[0,+\infty[$)
  ou $\RR$.
  S'il n'y a pas de précision, le terme {\it géodésique} désigne une droite géodésique. 
  On appellera {\it germe de rayon géodésique}, ou plus simplement {\it germe},  
une classe d'équivalence de
rayons géodésiques locaux pour la relation 
d'équivalence $r_1\sim r_2$ si $r_1$ et $r_2$ coïncident
sur un segment initial non vide et non réduit à un point. De même, la relation $r\sim r'$ si $r\cap r'$ est un rayon
géodésique local      
    est  une relation d'équivalence sur l'ensemble des rayons d'une géodésique locale.
On appelle {\it bout} (au sens de Freudhental) d'une géodésique locale une telle classe d'équivalence.    
    Une géodésique locale a deux bouts. 

Soient $(X,d)$ un espace métrique localement $\CAT(0)$ et $p:(\widetilde{X},\widetilde{d})\to(X,d)$ un revêtement universel localement isométrique. Si l'espace $(X,d)$
est complet (par exemple si $X$ est compact),  d'après le théorème de Cartan-Hadamard 
(voir \cite[Chap.~II.4]{BriHae99}), l'espace $(\widetilde{X},\widetilde{d})$ est complet et  $\CAT(0)$. Dans la suite, on supposera toujours que $(X,d)$ est 
complet. Nous noterons $\partial_\infty\widetilde{X}$ le bord à l'infini de $(\widetilde{X},\widetilde{d})$ et
$\partial^2_\infty\widetilde{X}=\partial_\infty\widetilde{X}\times\partial_\infty\widetilde{X}-\Delta$ (où $\Delta=\{(x,x),x\in\partial_\infty\widetilde{X}\}$), 
 muni de la topologie induite par la topologie produit. 

Soit $\ell:\RR\to (X,d)$ une géodésique locale. On dit que $\ell$ est \textit{périodique} s'il existe $T_0>0$ tel que pour 
tout $t\in\RR$, $\ell(t+T_0)=\ell(t)$, qu'elle est \textit{positivement (resp. négativement) périodique} s'il existe des réels 
strictements positifs $T$ et $T_0$ tels que $\ell(t+T_0)=\ell(t)$ pour tout $t>T$ (resp. $\ell(t-T_0)=\ell(t)$ pour tout $t<-T$). 
On dit que $\ell$ est 
     \textit{positivement (resp. négativement) compacte} si $\ell([0,+\infty[)$ (resp. $\ell(]-\infty,0])$) est 
compacte. Si $\ell$ est périodique ou positivement périodique, elle est aussi d'image compacte ou positivement compacte, mais la réciproque est fausse 
(voir le lemme \ref{grapheplongé}).

On note $\G$ l'ensemble des géodésiques locales paramétrées de $(X,d)$ et $\widetilde{\G}$ l'ensemble des géodésiques paramétrées de 
$(\widetilde{X},\widetilde{d})$ (puisque ($\widetilde{X},\widetilde{d})$ est $\CAT(0)$, les géodésiques locales de $(\widetilde{X},\widetilde{d})$ sont des géodésiques).
Alors $\RR$ agit sur $\G$ et 
$\widetilde{\G}$ par translations à la source. On note $\G_{np}$ et $\widetilde{\G}_{np}$ les quotients de $\G$ et $\widetilde{\G}$ pour ces actions, c'est-à-dire les
ensembles des 
géodésiques locales 
définies à changements d'origine près. Si $\ell$ appartient à $\G$, on note toujours $\ell$  sa classe dans $\G_{np}$. On appellera indistinctement 
{\it topologie des géodésiques} 
   la topologie compacte-ouverte sur $\G$ et $\widetilde{\G}$ et la topologie quotient de la topologie compacte-ouverte pour l'action par translations de $\RR$ sur 
   $\G_{np}$ et $\widetilde{\G}_{np}$.

 Si $\widetilde{\ell}$ appartient à $\widetilde{\G}$, elle est propre et converge vers deux points distincts de $\partial_\infty\widetilde{X}$ lorsque $t$ 
tend vers $\pm \infty$      
  et on note $E(\widetilde{\ell})=(\widetilde{\ell}(-\infty),\widetilde{\ell}(+\infty))\in\partial_\infty^2\widetilde{X}$.
  Si $F$ est un ensemble de géodésiques de $\widetilde{\G}$, on note 
$E(F)=\{E(\widetilde{\ell}),\widetilde{\ell}\in F\}$. 

 Enfin, si $\ell$ est une géodésique locale, on note $\ell^-$ la géodésique locale d'orientation opposée définie par $\ell^-(t)=\ell(-t)$. Si $F$ est un ensemble de 
 géodésiques locales, on appelle {\it support} de $F$ et on note $\operatorname{Supp}F$  la réunion des images des éléments de $F$.
 Puisqu'elles ne dépendent pas du paramétrages, ces définitions se prolongent aux géodésiques locales non paramétrées. 
On utilisera indistinctement les conventions précédentes pour tous les espaces métriques considérés.

%

\subsection{Ordres cycliques.}\label{ordrecyclique}

Soient $X$ un ensemble et $X^3$ l'ensemble des triplets de points de $X$. Un \textit{ordre cyclique (total)} sur $X$ est 
une fonction $o : X^3\to\{-1,0,1\}$ telle que : 

\medskip
\noindent
$\bullet$~ pour tout $(x,y,z)\in X^3$, $o(x,y,z)=0$ si et seulement si $\card(\{x,y,z\})\leqslant 2$;

\noindent
$\bullet$~ pour tout $(x,y,z)\in X^3,\, o(x,y,z)=o(y,z,x)=-o(x,z,y)$;

\noindent
$\bullet$~ pour tous les quadruplets $(x,y,z,t)$ d'\'el\'ements deux \`a deux distincts de $X$, si $o(x,y,z)=1$ et $o(x,z,t)=1$, alors $o(x,y,t)=1$. 

\medskip

L'ordre cyclique $\overline{o}$ \textit{opposé} à $o$ est défini par $\overline{o}(x,y,z)=-o(x,y,z)$ 
pour tout triplet $(x,y,z)\in X^3$. Si $X$ est un espace topologique et si $X^3$ est 
muni de la topologie produit, on demande que la restriction de $o$ au sous-ensemble des triplets d'éléments deux à deux distincts soit continue.
Si $X$ est muni 
d'une action d'un groupe $\Ga$, on demande que $o$ soit invariante pour l'action diagonale de $\Ga$ sur $X^3$. 
Nous renvoyons par exemple à \cite[2.3.1]{Wolf11} pour des compléments sur les ordres cycliques.

%
%

\subsection{Laminations géodésiques.}

Soient $(X,d)$ un espace métrique connexe localement $\CAT(0)$ (complet) et $p:(\widetilde{X},\widetilde{d})\to(X,d)$ un revêtement universel localement isométrique 
(dont le choix est indifférent).
On suppose qu'il
existe un ordre cyclique $o$ sur $\ddX$ (invariant par l'action diagonale
du groupe de revêtement de $p$ et dont la restriction aux triplets de points deux à deux distincts est continue). On dit alors que le triplet $(X,d,o)$ est un espace
métrique localement $\CAT(0)$ {\it enrubanné}. 
Deux couples de points $(x_1,y_1)$ et $(x_2,y_2)$ de $\ddX$ sont dits {\it entrelacés} si les points sont deux à deux distincts
et si $o(x_1,x_2,y_1)=-o(x_1,y_2,y_1)$. 

\bdefi\label{definonentrelacée} Deux géodésiques de $(\widetilde{X},\widetilde{d})$ sont dites entrelacées
si leurs couples de points à l'infini  sont entrelacés. Deux géodésiques locales  $\ell_1$ et $\ell_2$ de $(X,d)$ 
sont dites entrelacées s'il existe des relevés de $\ell_1$ et $\ell_2$ dans
$\widetilde{X}$ qui sont entrelacés. Si $\ell_1$ est entrelacée avec elle-même, elle est dite auto-entrelacée.
\edefi

Nous verrons dans la partie \ref{partienonentrelacée} des caractérisations des géodésiques locales non entrelacées ou auto-entrelacées, dans le cas particulier des 
surfaces munies de structures de demi-translation.

\bdefi\label{laminationplate1}
Une lamination géodésique (ou plus simplement lamination) de $(X,d,o)$ est 
un ensemble non vide $\Lambda$ de géodésiques locales de $(X,d)$, définies à changements d'origine près, dont les éléments sont appellés feuilles, tel que: 

\medskip

\noindent
$\bullet$~ aucune feuille de $\Lambda$ n'est auto-entrelacée;

\noindent
$\bullet$~ les feuilles de $\Lambda$ sont deux à deux non entrelacées;

\noindent
$\bullet$~ si $\ell$ appartient à $\Lambda$ alors $\ell^-$ aussi;

\noindent
$\bullet$~ $\Lambda$ est fermé pour la topologie des géodésiques.

\edefi

\rem Cette définition généralise les laminations hyperboliques (voir la partie \ref{laminationhyperbolique}). D'autres types de laminations, dans des arbres réels,
ont été introduites par M.~Bestvina, M.~Feighn et M.~Handel (voir \cite{BesFeiHan97}) et par T.~Coulbois, A.~Hilion et M.~Lustig (voir \cite{CouHilLus07}), mais ce ne sont 
pas des arbres enrubannés. Je remercie le rapporteur d'avoir suggéré d'étudier les liens entre les laminations telles qu'introduites dans \cite{CouHilLus07} et
cette définition.

\medskip

Soient $p':(X',d')\to(X,d)$ un revêtement localement isométrique de $(X,d)$, $\Lambda$ une lamination  de $(X,d)$ et $\Lambda'$ l'ensemble 
des relevés des feuilles de $\Lambda$ dans $X'$.   
   Alors il existe un revêtement universel $q:(\widetilde{X},\widetilde{d})\to(X',d')$ tel que $p=p'\circ q$, et les ensembles de relevés
des géodésiques locales de $\Lambda$ et $\Lambda'$ dans $\widetilde{X}$ sont égaux. Donc les géodésiques locales de
$\Lambda'$ sont deux à deux non entrelacées et aucune n'est auto-entrelacée. De plus, puisque $p'$ est localement isométrique, l'ensemble $\Lambda'$ est fermé, comme 
$\Lambda$, pour la topologie des géodésiques, et si $\ell'$ appartient à $\Lambda'$, alors $\ell'^-$ aussi. Donc $\Lambda'$ est une 
lamination   de $(X',d')$, que l'on appelle {\it image réciproque de $\Lambda$ par $p'$}. Dans cet article, l'espace $(\widetilde{X},\widetilde{d})$ sera généralement 
$\delta$-hyperbolique (au sens de Gromov) pour un certain $\delta\geqslant 0$
(voir \cite[Déf.~1.1~p.~399]{BriHae99}), et propre. De ce fait, pour tout couple de points distincts $(x,y)\in\partial_\infty^2\widetilde{X}$,
il existe au moins une géodésique de $(\widetilde{X},\widetilde{d})$
dont le couple de points à 
l'infini est $(x,y)$ (voir \cite[Lem.~3.2~p.~428]{BriHae99}). Dans la suite, on utilisera ce fait sans le mentionner. 
On rappelle un résultat bien connu, découlant du théorème d'Ascoli.

\blemm\label{compac}
Soit $(Z,d)$ un espace métrique propre. L'ensemble des g\'eod\'esiques de $(Z,d)$ dont l'origine appartient \`a un 
compact de $Z$ est compact pour la topologie compacte-ouverte. \cqfd
\elemm

\bcoro\label{Ffermé} Soit $F$ un ensemble de géodésiques définies à changements d'origines près, qui est fermé pour la topologie des géodésiques. 
Alors le support de $F$ est fermé dans $Z$.\cqfd
\ecoro

Du lemme \ref{compac} on déduit les trois lemmes élémentaires suivants, avec $(Z,d)$ un espace métrique géodésique, $\CAT(0)$ et propre.


\blemm\label{Xfermé}
Soit $Y$ un ensemble fermé de $\partial_\infty^2{Z}$ et $F$ l'ensemble de toutes les géodésiques $\ell$ de $(Z,d)$ telles que $E(\ell)$ appartient à $Y$.
Alors $F$ est fermé pour la topologie des géodésiques et donc $\operatorname{Supp}F$ est fermé dans $Z$.\cqfd
\elemm


\blemm\label{convergenceinfini}
Supposons de plus que l'espace $(Z,d)$ est $\delta$-hyperbolique, avec $\delta\geqslant 0$. Soient $({\ell}_n)_{n\in\NN}$ une suite de géodésiques 
(non paramétrées mais orientées) de $(Z,d)$ et pour tout $n\in\NN$, $(x_n,y_n)$ le couple de ses points à l'infini. 
Si la suite $(x_n,y_n)_{n\in\NN}$
converge vers $(x,y)$ dans $\partial_\infty^2{Z}$, alors quitte à extraire, la suite $({\ell}_n)_{n\in\NN}$ converge vers une géodésique de $(Z,d)$, 
pour la topologie des géodésiques, dont le couple de points à l'infini est $(x,y)$. \cqfd
\elemm

\dem Pour tout $n\in\NN$, on note $r_{x_n}=[A,x_n[$ et $r_{y_n}=[A,y_n[$, avec $A\in Z$ fixé. Alors les suites $(r_{x_n})_{n\in\NN}$ et $(r_{y_n})_{n\in\NN}$ convergent
vers les
  rayons $r_x=[A,x[$ et $r_y=[A,y[$ pour la topologie compacte-ouverte. Puisque $x\not=y$, il existe $T>0$ tel que $d(r_x(T),r_y(T))\geqslant 2\delta+3$. Or,
pour $n$ assez grand, on a  $\max\{d(r_{x_n}(T),r_{x}(T)),\\d(r_{y_n}(T),r_{y}(T))\}<1$, donc $d(r_{x_n}(T),r_{y_n}(T))\geqslant2\delta+1$. 
Puisque les triangles idéaux de côtés $r_{x_n}([0,+\infty[)$, $r_{y_n}([0,+\infty[)$ et $\ell_n(\RR)$ sont $\delta$-fins, pour tout $n\in\NN$, il existe un point $z_n$ de
$\ell_n(\RR)$ tel que $d(r_{x_n}(T),z_n)\leqslant \delta$, et donc $z_n$ appartient à la boule compacte de centre $r_x(T)$ et de rayon $\delta+1$. On
peut changer les origines
pour que $\ell_n(0)=z_n$ pour tout $n$ assez grand. Le lemme \ref{convergenceinfini} est alors une conséquence du lemme \ref{compac}.\cqfd

\blemm\label{caractérisation}
On suppose toujours que $(Z,d)$ est $\delta$-hyperbolique (avec $\delta\geqslant 0$). Soient $F$ un ensemble de géodésiques de $(Z,d)$, fermé pour la topologie des 
géodésiques,
et $Y_F$ l'ensemble des couples de points à l'infini des éléments de $F$. Alors $Y_F$ est fermé dans $\partial_\infty^2Z$.\cqfd
\elemm

\dem Soit $(x_n,y_n)_{n\in\NN}$ une suite  de $Y_F$ qui converge vers $(x,y)$
dans $\partial_\infty^2Z$. Pour tout $n\in\NN$, on note 
${\ell}_n$ un élément de $F$ dont le couple de points à l'infini est $(x_n,y_n)$. D'après le lemme \ref{convergenceinfini}, quitte à extraire,
la suite $(\ell_n)_{n\in\NN}$ converge, pour la topologie des géodésiques, vers une géodésique dont le couple de points à l'infini est $(x,y)$, et elle appartient à 
$F$ car $F$ est fermé pour cette topologie. Donc $Y_F$ est fermé. \cqfd

\subsection{Structures de demi-translation.}\label{structuresplates}

Soit $\Sigma$ une surface connexe, orientable, à bord (éventuellement vide). 
      Supposons que $\Sigma$ soit munie d'une métrique euclidienne sur $\Sigma-Z$, où $Z$ est une partie discrète de $\Sigma$.
Si l'holonomie de tout lacet fermé, $\C^1$ par morceaux, de $\Sigma-Z$ est contenue dans $\{\pm\Id\}$, on dit que deux vecteurs tangents $v_1$ et $v_2$ à $\Sigma$ sont 
{\it de même direction},  si $v_2$ est l'image de $\pm v_1$ par holonomie le long d'un chemin $\C^1$ par morceaux de $\Sigma-Z$ entre les points bases de $v_1$ et $v_2$. 
Cette définition ne dépend pas du choix du chemin. On dit qu'un chemin, ou une réunion de chemins, $\C^1$ par morceaux est {\it de direction constante},
si tous ses vecteurs tangents, aux points appartenant à $\Sigma-Z$, sont de même direction.

\bdefi
Une structure de demi-translation (ou structure plate à holonomie $\{\pm\Id\}$) sur une surface $\Sigma$ est la donn\'ee d'une partie discrète $Z$ de $\Sigma$
(\'eventuellement vide) et d'une m\'etrique euclidienne sur $\Sigma-Z$ \`a singularit\'e conique d'angle $k_z\pi$ en $z\in Z$, avec $k_z\in\NN$ et $k_z\geqslant 3$
  si
$z\in Z-Z\cap\partial\Sigma$ et $k_z\geqslant 2$ si $z\in Z\cap\partial\Sigma$, telle que l'holonomie de tout lacet fermé, $\C^1$ par morceaux, 
de $\Sigma-Z$ soit contenue dans 
$\{\pm\Id\}$ et telle que la réunion des composantes connexes du bord soit de direction constante.
\edefi

Nous donnerons dans la partie \ref{diffquadholo} des caractérisations des structures de demi-translation. Une structure de demi-translation définit une
distance géodésique qui est localement $\CAT(0)$,
que l'on notera $d$. On supposera toujours que $(\Sigma,d)$ est complet.
%
%
Nous appellerons {\it géodésiques locales plates} les géodésiques locales d'une structure de demi-translation.
Une application continue $\ell:\RR\to\Sigma$ est une géodésique locale plate (pour $d$) si
et seulement si elle vérifie (voir \cite[Th.~5.4~p.24]{Strebel84} et \cite[Th.~8.1~p.~35]{Strebel84}):

\medskip
\noindent
$\bullet$~ pour tout $t\in\RR$ tel que $\ell(t)$ n'appartient pas à $Z$, il existe un voisinage $V$ de $t$ dans $\RR$ tel que $\ell_{|V}$ soit un segment euclidien
 (donc de direction constante).
 
\noindent
$\bullet$~ pour tout $t\in\RR$, si $\ell(t)$ appartient à $Z$, alors les
 deux angles définis par les germes $\ell([t,t+\varepsilon[)$ et $\ell(]t-\varepsilon,t])$, avec $\varepsilon>0$ assez petit, mesurés dans chaque composante connexe
 de $U-\ell(]t-\varepsilon,t+\varepsilon[)$, où $U$ est un voisinage assez petit de $\ell(t)$, sont supérieurs ou égaux à $\pi$.
\begin{center}
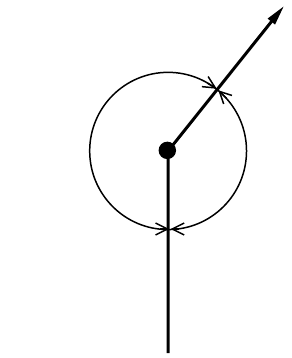
\end{center}

\medskip

Une \textit{liaison de singularités} d'une structure de demi-translation est l'image d'un segment
géodésique local entre deux singularités (éventuellement confondues) qui ne rencontre pas de 
singularité en dehors de ses extrémités. Puisque les géodésiques locales ne peuvent changer de direction qu'aux singularités, une liaison de singularités est de
direction constante.

Soit $p:(\widetilde{\Sigma},\widetilde{d})\to(\Sigma,d)$ un revêtement universel localement isométrique de $(\Sigma,d)$. 
L'espace métrique $(\widetilde{\Sigma},\widetilde{d})$ est une surface munie d'une structure de demi-translation qui est $\CAT(0)$ et propre.
Nous noterons $\partial_\infty\widetilde{\Sigma}$ son bord à l'infini.

\brema\label{ordresurface}{\rm Si $(\Sigma,m)$ est une surface connexe, orientable, sans bord, munie d'une métrique hyperbolique complète, et si $p:(\widetilde{\Sigma},
\widetilde{m})
\to(\Sigma,m)$ est un revêtement universel riemannien, alors le bord à l'infini de $\widetilde{\Sigma}$ est homéomorphe à  $\SSS^1$
(voir \cite[Ex.~8.11~(2)~p.265]{BriHae99}). Donc les ordres opposés usuels de $\SSS^1$ induisent deux ordres opposés sur $\ddp$. De même, si le bord de
$\Sigma$ n'est pas vide et est totalement
géodésique, alors $(\Sigma,m)$ se plonge isométriquement dans une (unique à isométrie près) surface orientable et hyperbolique (complète), sans bord, soit $(\Sigma',m')$, 
qui se rétracte par déformation forte sur $\Sigma$. Alors $\partial_\infty\widetilde{\Sigma}$ est un sous-ensemble fermé du bord à l'infini
$\partial_\infty\widetilde{\Sigma}'$ d'un revêtement universel (riemannien) de $\Sigma'$ (voir \cite[Ex.~8.11~(4)~p.266]{BriHae99}) et les ordres
opposés
sur $\partial_\infty\widetilde{\Sigma}'$ induisent par restrictions des ordres opposés sur $\partial_\infty\widetilde{\Sigma}$. 
  Si de plus $\Sigma$ est compacte et $\chi(\Sigma)<0$, et si on note $\Gamma_\Sigma$ le groupe de revêtement de $p:\widetilde{\Sigma}\to\Sigma$, alors  
  il existe un unique homéomorphisme $\Gamma_\Sigma$-équivariant entre $\partial_\infty\Gamma_\Sigma$ et le bord à l'infini de $(\widetilde{\Sigma},\widetilde{d})$ 
  pour n'importe quelle distance localement $\CAT(0)$ $d$ sur $\Sigma$, où $p:(\widetilde{\Sigma},\widetilde{d})\to(\Sigma,d)$ est localement isométrique, par 
  lequel on les identifie (voir par exemple \cite[§2]{Bonahon91}). 
  On peut donc définir cet ordre cyclique sur le bord à l'infini de $(\widetilde{\Sigma},\widetilde{d})$ pour n'importe quelle distance géodésique
  localement $\CAT(0)$ complète $d$ sur
  $\Sigma$. Les espaces métriques $(\widetilde{\Sigma},\widetilde{d})$ et $\revetm$ sont de plus quasi-isométriques, donc il existe $\delta>0$ tel que
  $(\widetilde{\Sigma},\widetilde{d})$ est $\delta$-hyperbolique. 


Dans la suite, on ne considérera que l'ordre cyclique $o_\infty$ ainsi défini sur le bord à l'infini d'un revêtement universel d'une surface compacte
munie d'une structure de demi-translation (et son opposé), et on ne 
le précisera pas. Si on choisit une orientation de $\widetilde{\Sigma}$ et si $or_x$ désigne l'ordre cyclique induit par l'orientation de $\widetilde{\Sigma}$ sur les
germes de rayons géodésiques issus d'un point $x$ de $\widetilde{\Sigma}$, quitte à remplacer l'ordre cyclique $o_\infty$ par son opposé, 
  si $a,b,c$ sont trois points deux à deux distincts de $\partial_\infty\widetilde{\Sigma}$
et si les germes $r_a,r_b,r_c$ des rayons géodésiques $[x,a[,[x,b[,[x,c[$ sont deux à deux distincts, on a}
\begin{equation}\label{eq:1}
or_x(r_a,r_b,r_c)=o_\infty(a,b,c).
\end{equation}
\erema

%

\subsection{Différentielles quadratiques holomorphes et structures plates dirigées.}\label{diffquadholo}

Soit $\Sigma$ une surface connexe, orientable, à bord (éventuellement vide). Rappelons qu'une structure de surface de Riemann à bord sur $\Sigma$ est la donnée d'un
atlas maximal de carte 
$(z_\nu : U_\nu \to V_\nu)_{\nu\in\tau}$, où $\{U_\nu\}_{\nu\in\tau}$ est un recouvrement ouvert de $\Sigma$ et  pour tout $\nu\in\tau$, l'application 
$z_\nu$ est un homéomorphisme sur un ouvert (relatif) $V_\nu$ de 
$\{z\in\CC\,:\,\Im(z)\geqslant 0\}$, tel que les changements de cartes soient holomorphes (voir par exemple \cite[Déf.~1.2~p.~2]{Strebel84}). Le bord de $\Sigma$ est l'ensemble des points  $x$ de $\Sigma$ tels qu'il existe une carte de cette structure contenant $x$, envoyant
$x$ sur l'axe réel, ce qui ne dépend pas d'une telle carte. Muni de la restriction  des cartes, le bord de $\Sigma$ est une variété analytique réelle de dimension $1$, non connexe a priori.


\bdefi\label{defdemitranslation}
Une différentielle quadratique holomorphe sur ${\Sigma}$ est la donnée d'une structure de surface de Riemann à bord 
$X$ sur ${\Sigma}$ et d'une famille $q$ de fonctions holomorphes $(q_{U,z}:U\to\CC)_{(U,z)\in X}$ telles que pour toutes les cartes locales $(U,z)$ et $(V,w)$ de $X$,   
on a

\medskip
\noindent
$(1)$~  $q_{U,z}dz^2=q_{V,w}dw^2$ sur $U\cap V$; 

\noindent
$(2)$~ l'image $q_{U,z}(U\cap\partial\Sigma)$ de
 $U\cap\partial\Sigma$ est contenue dans $\RR^-=\;]-\infty,0]$.
\edefi

\rem\label{remarquedefdemitranslation} D'après la condition $(1)$, il suffit de définir les différentielles quadratiques holomorphes $q_{U,z}dz^2$
pour une sous-famille de cartes locales qui 
recouvre $\Sigma$. Cette famille peut être finie si $\Sigma$ est compacte. Si $(U,z)$ et $(V,w)$ sont des cartes locales et si $y\in U\cap V$, alors 
$q_{V,w}(y)=q_{U,z}(y)\frac{dz(v)^2}{dw(v)^2}$. Donc $q_{V,w}(y)\in\RR^-$ si et seulement si $q_{U,z}(y)\in\RR^-$, et il suffit
de vérifier la deuxième condition pour une sous-famille de cartes locales qui recouvre le bord. Cette condition impose que l'image d'un segment du bord par un 
paramètre naturel (voir par exemple \cite[5.1~p.~21]{Strebel84}), 
défini par $q$ au voisinage d'un point régulier, soit vertical.

\medskip

On note $\widetilde{\Q}(\Sigma)$ l'ensemble des diff\'erentielles quadratiques holomorphes sur $\Sigma$ et seulement $q$ une différentielle quadratique holomorphe
$(X,q)$ lorsque la structure 
sous-jacente est sous-entendue.

\medskip

\rem Soient $p:\Sigma'\to\Sigma$ un rev\^etement de $\Sigma$ et $(X,q)$ une différentielle quadratique holomorphe sur ${\Sigma}$. Le couple de la structure de surface
de Riemann à bord $X'$ tirée en arrière de $X$ sur $\Sigma'$ et de la famille $q'=(q_{p(V),z\circ(p_{|V})^{-1}}\circ p_{|V})_{(V,z)}$, où $(V,z)$
parcourt l'ensemble des cartes locales de $X'$ telles que $p_{\mid V}$ est un difféomorphisme sur son image définit une différentielle quadratique holomorphe sur 
$\Sigma'$, dite
relevée de $(X,q)$, d'après la remarque précédente.

\medskip

Le groupe $\Diff_0({\Sigma})$ des isotopies de ${\Sigma}$ qui fixent point par point $\partial{\Sigma}$ et sont homotopes à l'identité par des applications
qui fixent point par point $\partial\Sigma$, agit sur $\widetilde{\Q}(\Sigma)$ par pr\'ecomposition des cartes 
et des expressions des diff\'erentielles quadratiques holomorphes dans ces cartes. On note $\Q(\Sigma)$ le quotient de $\widetilde{\Q}(\Sigma)$ par cette action.

%
%
%
%

\bdefi
Une structure plate dirigée sur $\Sigma$ est la donn\'ee:

\medskip

\noindent
$\bullet$~ d'une partie discrète $Z$ de $\Sigma$ (\'eventuellement vide);

\noindent
$\bullet$~ d'une m\'etrique euclidienne \`a singularit\'e conique d'angle $k_z\pi$ en $z\in Z$, avec $k_z\in\NN$ et $k_z\geqslant 3$
  si
$z\in Z-Z\cap\partial\Sigma$ et $k_z\geqslant 2$ si $z\in Z\cap\partial\Sigma$;

\noindent
$\bullet$~ d'un champ de droites 
parall\`eles $V$ sur $\Sigma-Z$, tel que pour tout $x\in\partial\Sigma-Z\cap\partial\Sigma$, la droite $V(x)$ soit contenue dans $T_x\partial\Sigma$.

\edefi  

Le groupe $\Diff_0({\Sigma})$ agit sur l'ensemble des structures plates dirigées
sur $\Sigma$ par tirés en arrière de la métrique et du champ de droites parallèles. On note $\widetilde{\Plat}(\Sigma)$ le quotient de l'ensemble des structures
plates dirigées par cette action (voir la fin de cette partie pour une explication de la notation).

\medskip

La donn\'ee d'une différentielle quadratique holomorphe d\'efinit une m\'etrique 
euclidienne \`a singularit\'es coniques du bon type sur ${\Sigma}$ et, un champ de droites parallèles (verticales) sur le complémentaire des singularités.
L'application qui à une différentielle quadratique holomorphe associe cette structure plate dirigée est
$\Diff_0(\Sigma)$-équivariante et définit une application de $\Q(\Sigma)$ dans $\widetilde{\Plat}(\Sigma)$, qui est une bijection car $\Sigma$ est orientable, par laquelle
on identifie $\Q(\Sigma)$ et $\widetilde{\Plat}(\Sigma)$.



\medskip

Le groupe $(\RR,+)$ agit 
  sur $\Q(\Sigma)$ par l'application $(\theta,
(X,q))\mapsto(X,e^{2i\theta}q)$. Si $(X,q)\in\Q(\Sigma)$, on note $(X,[q])$, ou plus simplement $[q]$, sa classe d'équivalence pour cette action. De même, il 
agit sur l'ensemble $\widetilde{\Plat}(\Sigma)$ par rotation des champs de vecteurs parallèles. De plus, la bijection entre $\Q(\Sigma)$ et 
$\widetilde{\Plat}(\Sigma)$ est équivariante pour ces actions et définit une bijection entre les espaces quotients. 

Une structure plate dirigée est en particulier une structure de demi-translation, d'après l'existence d'un champ de droites parallèles. Il existe donc une 
application
d'oubli de l'ensemble des structures plates dirigées dans l'ensemble des structures de demi-translation, qui est invariante pour l'action par rotation de 
$\RR$.
Elle définit donc une application du quotient de $\widetilde{\Plat}(\Sigma)$ par cette action à valeurs  dans l'ensemble des structures de demi-translation,
qui est une bijection. 
On identifie ainsi l'ensemble des structures de demi-translation avec l'espace quotient $\Plat(\Sigma)=\RR\backslash\widetilde{\Plat}(\Sigma)$, et on notera 
$[q]$, avec $q\in\Q(\Sigma)$, une structure de demi-translation. On ne confondra pas une surface de demi-translation qui peut être définie comme un couple $(\Sigma,q)$ 
où $q\in\Q(\Sigma)$, et la structure de demi-translation $[q]$ qui lui est associée. 
%
%
%

\section{Comportement des géodésiques locales d'une structure de demi-translation.}
Soient $\srfce$ une surface connexe, orientable, à bord (éventuellement vide), munie d'une      
     structure de demi-translation et $p:\revet\to\srfce$ un revêtement universel localement isométrique. On note $or_{x}$ l'ordre cyclique sur
l'ensemble des germes de rayons
issus d'un point $x$ induit par un choix de l'orientation de la surface et $o_\infty$ l'ordre sur $\partial_\infty\widetilde{\Sigma}$ défini dans la remarque 
\ref{ordresurface}. 
%
%
%
%
%
%
%
%
\subsection{Comportement de deux géodésiques locales d'une structure de demi-translation.}\label{partienonentrelacée}

  Soient $\ell_1$ et ${\ell}_2$ deux géodésiques locales paramétrées de $\srfce$ telles qu'il existe $I\subset\RR$ un intervalle maximal compact, non vide mais éventuellement
réduit à un point, tel que 
$\ell_{1|I}=\ell_{2|I}$. Dans la suite, on notera alors $t_1\leqslant t_2$ tels que $I=[t_1,t_2]$, $x_1=\ell_1(t_1)$, $x_2=\ell_1(t_2)$, et $r_1^-,r_2^-,r_1^+,r_2^+$ les germes des rayons
géodésiques
${\ell}_1(]t_1-\varepsilon,t_1])$, ${\ell}_2(]t_1-\varepsilon,t_1])$, ${\ell}_1([t_2,t_2+\varepsilon[)$, 
${\ell}_2([t_2,t_2+\varepsilon[)$ issus  de $x_1$ et 
$x_2$, et si $t_1\not=t_2$, $r_0^-,r_0^+$ les germes de
${\ell}_1([t_1,t_1+\varepsilon[)={\ell}_2([t_1,t_1+\varepsilon[)$ et ${\ell}_1(]t_2-\varepsilon,t_2])=
{\ell}_2(]t_2-\varepsilon,
t_2])$, pour un $\varepsilon>0$ assez petit (voir la figure du lemme).

\medskip

Soient $\widetilde{\ell}_1$ et $\widetilde{\ell}_2$ deux géodésiques de $\revet$. Par définition des ordres cycliques, $\widetilde{\ell}_1$ et $\widetilde{\ell}_2$
sont non entrelacées si et seulement si $\widetilde{\ell}_1^i$ et $\widetilde{\ell}_2^j$    
    sont non entrelacées (pour tout $i,j\in\{+,-\}$). Donc si leurs images sont confondues le long d'un 
segment ou un rayon géodésique, on peut les paramétrer pour que leurs restrictions au segment ou rayon soient égales.
Or, puisque l'espace $\revet$ est $\CAT(0)$ et les géodésiques sont continues, si l'intersection des images de $\widetilde{\ell}_1$ et $\widetilde{\ell}_2$
n'est pas vide, c'est un point isolé ou une réunion connexe et fermée de liaisons de singularités et (éventuellement) d'images de rayons géodésiques
 d'origine singulière. De plus, les géodésiques ne peuvent  changer de direction qu'aux singularités, donc si les
géodésiques n'ont pas la même image, les (ou l'unique) extrémités de $\widetilde{\ell}_1(\RR)\cap\widetilde{\ell}_2(\RR)$ sont des singularités. D'après la formule 
\eqref{eq:1}, on a le lemme suivant et son corollaire.

\blemm\label{nonentrelacéeordre}
Deux géodésiques $\widetilde{\ell}_1$ et $\widetilde{\ell}_2$ de $\revet$ sont non entrelacées si et seulement si, en reprenant les notations
précédentes, l'une des conditions exclusives suivantes est vérifiée:

\medskip

\noindent
$(1)$~ Leurs images sont disjointes;

\noindent
$(2)$~ Leurs images sont confondues au moins le long d'un rayon géodésique;

\noindent
$(3)$~ Leurs images s'intersectent en un unique point $x$ et $or_{x}(r_1^+,r_2^+,r_1^-)=or_{x}(r_1^+,r_2^-,r_1^-)$;

\noindent
$(4)$~ quitte à changer d'orientation et d'origines, leurs restrictions à un segment géodésique maximal $[x_1,x_2]$ non réduit à un point sont égales, 
et $or_{x_2}(r_1^+,r_2^+,r_0^+)=-or_{x_1}(r_1^-,r_2^-,r_0^-)$.
\begin{center}
 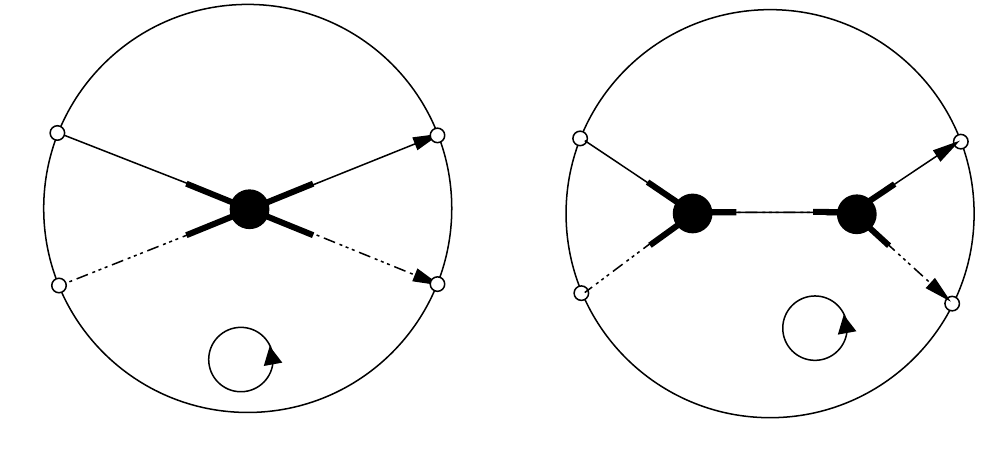
\end{center}
\elemm

\dem Le complémentaire de $\widetilde{\ell}_1(\RR)$ dans $\widetilde{\Sigma}$ a deux composantes connexes, et $\widetilde{\ell}_1$ et $\widetilde{\ell}_2$ sont non
entrelacées si et seulement si $\widetilde{\ell}_2(\RR)$ est contenue dans l'adhérence d'une seule. Une analyse de cas montre que cela correspond aux conditions du lemme.\cqfd

\medskip

Si $\ell_1$ et $\ell_2$ sont des géodésiques locales d'images disjointes de $\srfce$, les images de leurs relevés sont deux à deux disjointes et d'après le 
lemme \ref{nonentrelacéeordre}, elles sont non entrelacées. Si elles s'intersectent, puisque $p$ préserve l'orientation, on a les deux corollaires suivants.

\bcoro\label{coro1}
Si $\ell_1$ et $\ell_2$ s'intersectent, elles sont non entrelacées si et seulement si 
l'intersection des images de $\ell_1$ et $\ell_2$ est
une réunion fermée non vide de singularités, de liaisons de singularités et (éventuellement) d'images de rayons géodésiques locaux d'origine singulière 
(qui n'est pas nécessairement connexe), 
et si de plus elles vérifient: pour chaque intervalle maximal $I\subset\RR$ tel que, quitte à changer d'orientation et d'origines, on a $\ell_{1|I}=\ell_{2|I}$,
en reprenant 
les notations précédentes, l'une des conditions exclusives suivantes est vérifiée.

\medskip
\noindent
$\bullet$~ $I$ est de la forme $[t,+\infty[$ ou $]-\infty,t]$, avec $t\in\RR$;

\noindent
$\bullet$~ $I=\{t\}$, et si $x=\ell_1(t)$, on a $or_{x}(r_1^+,r_2^+,r_1^-)=or_{x}(r_1^+,r_2^-,r_1^-)$;

\noindent
$\bullet$~ $I=[t_1,t_2]$ et si $x_1=\ell_1(t_1)$ et $x_2=\ell_1(t_2)$, on a $or_{x_2}(r_1^+,r_2^+,r_0^+)=-or_{x_1}(r_1^-,r_2^-,r_0^-)$.\cqfd
\ecoro

\bcoro\label{coro2}
Si $\ell_1$ et $\ell_2$ sont deux géodésiques locales de $\srfce$ qui s'intersectent en un point $x$ qui n'est pas une singularité, 
et dont les images ne sont pas confondues au voisinage de $x$, alors $\ell_1$ et $\ell_2$ sont entrelacées.\cqfd
\ecoro
%
%
%

\subsubsection{Géodésiques locales non auto-entrelacées.}
On suppose dans cette partie que la surface $\Sigma$ est compacte de genre $g\in\NN$, et on note $b\in\NN$ le nombre de composantes connexes de son bord.
Soient $x_0$ un point de $\Sigma$ et  $B=(c_1,\dots,c_n)$ une suite finie
d'arcs orientés simples fermés en $x_0$, que l'on peut supposer paramétrés par $[0,1]$ avec $c_i(0)=c_i(1)=x_0$ pour tout $i=1,\dots,n$, 
tels que pour tout $k\not =p$, on ait $c_k([0,1])\cap c_p([0,1])=\{x_0\}$. On note, par abus, encore $B$ la réunion $c_1([0,1])\cup\dots\cup c_n([0,1])$. 
On munit $\Sigma$ de la distance associée à une métrique riemannienne à bord totalement géodésique
quelconque sur $\Sigma$. 
Soient $F_1,\dots,F_k$ 
les complétés des composantes connexes de $\Sigma-B$ (pour la métrique riemannienne induite). Pour tout $i=1,\dots,k$, $F_i$ est une surface (topologique) à bord.

\blemm\label{Euler}
La caractéristique d'Euler de $\Sigma$ vérifie
$\chi(\Sigma)=1-n+\sum_1^k\chi(F_i)$.
\elemm

\dem Pour tout complété $F_i$ d'une composante connexe de $\Sigma-B$, si $c$ est une composante du bord de $F_i$, on note $V_{i,c}\simeq c\times[0,\varepsilon]
\subset F_i$,
avec $\varepsilon>0$ assez petit, un voisinage tubulaire de $c$ dans $F_i$, et $F_i'$ le complété de $F_i-\underset{c}\amalg V_{i,c}$ pour la métrique induite, où 
$c$ parcourt
l'ensemble des composantes du bord de $F_i$. Puisque $F_i $ se rétracte par déformation forte sur $F_i'$,
on a $\chi(F_i')=\chi(F_i)$. De plus, pour tout $i=1,\dots,k$, la sous-surface à bord $F_i'$ se plonge dans $\Sigma$. On note $U$ la réunion des images dans $\Sigma$ 
des
voisinages $V_{i,c}$, pour
$i=1,\dots,k$ et $c$ une composante du bord de $F_i$. Alors $U$ se rétracte par déformation forte sur $B$, donc $\chi(U)=\chi(B)$, et $\Sigma$ est la réunion de $U$ et des
images des $F_i'$ dans $\Sigma$ (que l'on note toujours $F_i'$ car elles se plongent dans $\Sigma$). Pour toute sous-surface à bord $F_i'$, l'intersection $F_i'\cap U$ 
est homéomorphe à un cercle, de caractéristique d'Euler nulle, et 
$\chi(\Sigma)=\chi(U)+\sum^k_1\chi(F_i')=\chi(B)+\sum_1^k\chi(F_i)$. Enfin 
$\chi(\Sigma)=1-n+\sum_1^k\chi(F_i)$. \cqfd

\medskip

On dit que la suite finie $(c_1,\dots,c_n)$ est un \textit{bouquet plongé essentiel} en $x_0$ si de plus aucun arc $c_i$ n'est librement homotope à un point, et
si pour $k,p$ distincts dans $\{1,...,n\}$,
$c_k$ n'est homotope ni à $c_p$ ni à $\overline{c_p}$ relativement à $x_0$ (où $\overline{c_p}$ est défini par $\overline{c_p}(t)=c_p(1-t)$).

\begin{center}
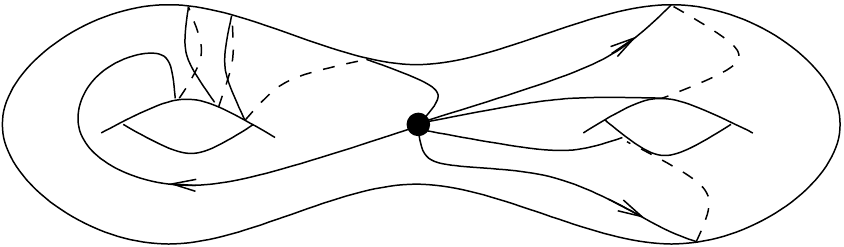
\end{center}

\blemm\label{Bouquet}
Si $B=(c_1,\dots,c_n)$ est un bouquet plongé essentiel en $x_0$, alors $n\leqslant 6g+3b-3$.
\elemm

\dem Soient $F_1,...,F_k$ comme ci-dessus. Pour tout $i=1,\dots,k$, si $F_i$ n'est pas un disque,
alors $\chi(F_i)\leqslant 0$. De plus, si $F_i$ est un disque, il est
bordé par au moins 3 arcs de la suite finie, car sinon il existerait deux arcs homotopes relativement à $x_0$, ou un arc homotope à un point. De plus, chaque arc 
ne borde,
au plus, que deux composantes connexes de $\Sigma-B$. Quitte à renuméroter, on peut supposer que les composantes $F_1,...,F_p$  (où $p\leqslant k$) sont tous les disques. Alors
$p\leqslant\frac{2n}{3}$. Or, d'après le lemme \ref{Euler}, on a $\chi(\Sigma)=1-n+\sum_1^k\chi(F_i)$. Puisque la caractéristique d'Euler d'un disque est $1$
et $\chi(F_i)\leqslant 0$ si $i>p$, on a
$\chi(\Sigma)\leqslant 1-n+\frac{2n}{3}=1-\frac{n}{3}$. Puisque $\chi(\Sigma)=2-2g-b$, on a bien $n\leqslant 6g+3b-3$.\cqfd

\bcoro\label{coroBouquet}
Si $x$ et $y$ sont deux points distincts de $\Sigma$ et si $(c_1,\dots,c_n)$ est une suite finie d'arcs simples joignant $x$ à $y$ qui sont deux à deux non homotopes relativement à $x$ et $y$, 
et sont deux à deux d'images disjointes en dehors de $x$ et $y$, alors $n\leqslant 6g+3b-2$.    
\ecoro

\dem En écrasant $c_1([0,1])$ en un point par une application continue de $\Sigma$ qui est un homéomorphisme en dehors de $c_1([0,1])$, on obtient un bouquet plongé
essentiel de $n-1$
arcs fermés simples en un point 
de $\Sigma$. D'après le lemme \ref{Bouquet}, on a $n-1\leqslant 6g+3b-3$, donc $n\leqslant 6g+3b-2$. \cqfd

\medskip

On considère de nouveau une structure de demi-translation $[q]$ sur $\Sigma$. Puisque $\Sigma$ est compacte,
le nombre $k$ de singularités de $[q]$
est fini. Si $\ell:\RR\to\Sigma$ est une géodésique locale et 
 $R$ est l'image d'un segment géodésique local de longueur finie, on dit que $\ell$ \textit{parcourt} $R$  
 s'il existe au moins un intervalle compact $I$ de $\RR$ tel que $\ell(I)=R$. Nous aurons besoin du lemme bien connu suivant.

\blemm\label{homotopie}
Soient $(X,d)$ un espace géodésique localement $\CAT(0)$ et $(x,y)$ un couple de points de $X$. Alors il existe un unique
 segment géodésique dans chaque classe d'homotopie relative d'arcs joignant $x$
à $y$. \cqfd
\elemm

%

\blemm\label{fini}
Soit $F$ une famille de géodésiques locales non auto-entrelacées de $\srfce$ qui sont deux à deux non entrelacées.
Alors il y a au plus $k^2\cdot (6g+3b-2)$ liaisons de singularités qui sont parcourues par une géodésique locale de $F$.
\elemm

%
%

\dem Soient $h$ et $h'$ deux liaisons de singularités distinctes joignant deux singularités $x$ et $y$ de $[q]$. Puisqu'une liaison de singularités est de direction
constante, si $h$ et $h'$ s'intersectent en dehors de leurs extrémités, cette intersection est transverse. Donc, d'après le corollaire \ref{coro2}, si $h$ et $h'$
sont parcourues par des géodésiques de $F$, elles ne s'intersectent pas en dehors de leurs extrémités. De plus, d'après le lemme \ref{homotopie}, elles ne sont pas
homotopes
relativement à $x$ et $y$. Donc, quitte à les ordonner, la famille des liaisons de singularités joignant $x$ et $y$ est une suite d'arcs simples satisfaisant
les hypothèses du corollaire \ref{coroBouquet}, donc il y en a au plus $6g+3b-2$.     

   Puisque $[q]$ n'a que $k$ singularités, 
il y a $k^2$ couples de singularités. Donc il y a au plus
$k^2\cdot (6g+3b-2)$ liaisons de singularités qui sont parcourues par une géodésique de $F$.
\cqfd

\blemm\label{directionconstante}
Soit $\ell:\RR\to\Sigma$ une géodésique locale non auto-entrelacée qui n'est pas positivement compacte. Il existe $T>0$ tel que $\ell_{|[T,+\infty[}$ 
est de direction constante.
\elemm

\dem
D'après le lemme \ref{fini} appliqué à $F=\{\ell\}$, il n'y a qu'un nombre fini de liaisons de singularités qui sont parcourues par $\ell$. Puisque toutes les liaisons de singularités sont compactes, la réunion 
des liaisons de singularités parcourues par $\ell$ est compacte.
Donc $\ell$ ne rencontre pas de singularité à partir d'un certain temps. Elle est alors de direction constante à partir de ce temps.\cqfd

%
%
%

\medskip

Si $q$ est une différentielle quadratique holomorphe sur $\Sigma$, une \textit{trajectoire} de $q$ est un segment, un rayon ou une droite
géodésique locale, verticale pour $q$, 
ne rencontrant pas de
singularité, maximale (voir \cite[Déf.~5.5.3~p.25]{Strebel84}).
Si $\ell:I\to\Sigma$ est
    une trajectoire ($I$ est alors un intervalle ouvert), un {\it rayon de trajectoire} de $\ell$ est 
une application $r$ de la forme $\ell_{|[T,a[}$ ou $\ell_{|]a,T]}$, avec $a\in\RR\cup\{\pm\infty\}$ et $T\in I$ tels que $[T,a[\;=I\cap[T,+\infty[$ ou $]a,T]=I\;\cap\;
]-\infty,T]$.
Cette définition ne dépend pas du paramétrage de $\ell$, quitte
à changer les valeurs de $T$ et $a$.
Alors si $a<+\infty$ et en reparamétrant $r$ sur $[0,|a-T|[$, on dit que $r$ {\it converge vers un zéro} $x$ de $q$
si $r(t)\underset{t\to |a-T|^-}\longrightarrow x$.
Si $a=\pm\infty$, on dit que $r$ est {\it récurrent} s'il est dense dans un {\it domaine}, c'est-à-dire l'adhérence d'un ouvert connexe de $\Sigma$ bordé par
les images de feuilles 
verticales périodiques rencontrant au moins une singularité, s'il n'est pas égal à $\Sigma$, du feuilletage vertical de $q$ (voir \cite[p.53]{Strebel84}). 
Un domaine n'est pas toujours une sous-surface de $\Sigma$ car son bord peut ne pas être une sous-variété de $\Sigma$.   
\blemm\label{structuredestrajectoires}(Voir \cite[11.4~p.53]{Strebel84})
Les trajectoires d'une différentielle quadratique holomorphe $q$ d'une surface compacte, à bord  (éventuellement vide) $\Sigma$ sont classées dans les familles suivantes:

\medskip
\noindent
$(1)$~  Les trajectoires périodiques. Elles feuillettent des cylindres plats ouverts.

\noindent
$(2)$~ Les trajectoires non périodiques.

\noindent
$(2.a)$~  Les trajectoires critiques, dont au moins un rayon converge vers une singularité. Il n'y a qu'un nombre fini de trajectoires critiques.
Si un rayon d'une trajectoire critique ne converge pas vers une singularité, il admet un sous-rayon récurrent. 

\noindent
$(2.b)$~  Les spirales, c'est-à-dire les trajectoires dont tous les rayons sont récurrents. Si $\alpha$ est une spirale,  
l'adhérence $A$ de l'image de $\alpha$ est un domaine, bordé par des réunions connexes d'images de trajectoires critiques de longueur finie et de singularités,
s'il n'est pas égal à $\Sigma$.    
\elemm

\medskip

\rm{Puisque ce sont les seules possibilités et que deux trajectoires de la liste ne peuvent avoir de point en commun à moins  qu'elles ne coïncident, 
la surface $\Sigma$ privée des zéros de $q$ et des liaisons de singularités verticales de $q$ se divise entre des cylindres plats ouverts (feuilletés
par des trajectoires périodiques) et les intérieurs des domaines qui sont les adhérences de chacune des spirales qu'ils contiennent.}

\bcoro\label{corostructuredestrajectoires}
Si $\ell$ est une géodésique locale de $\srfce$ telle que  $\ell_{|[T,+\infty[}$ ne rencontre pas de singularité, et est donc de direction constante, quitte à choisir un autre représentant de $[q]$, on peut supposer que
$\ell_{|[T,+\infty[}$ est de direction verticale pour $q$, et donc que c'est un rayon de trajectoire de $q$. D'après la classification ci-dessus, puisque le rayon 
$\ell_{|[T,+\infty[}$ ne rencontre pas de singularité, il est périodique ou récurrent et dans le second cas, son image est dense dans un domaine
 du feuilletage vertical de $q$.
\ecoro
\section{Image au bord à l'infini d'une famille de géodésiques.}\label{partiebordinfini}
%

\subsection{Rappels sur les laminations géodésiques hyperboliques.}\label{laminationhyperbolique}

Soient $\Sigma$ une surface connexe, orientable, à bord (éventuellement vide), et $m$ une métrique hyperbolique à bord
totalement géodésique (complète) sur $\Sigma$.
Une {\it lamination géodésique hyperbolique} (ou plus simplement {\it lamination hyperbolique}) $\Lambda_m$ de $\srfcem$ est une partie non vide fermée de $\Sigma$,
réunion d'images
de géodésiques locales hyperboliques, simples et deux à deux disjointes. L'espace $\srfcem$ est un espace métrique localement $\CAT(-1)$ et donc 
localement $\CAT(0)$. Soit $p:\revetm\to\srfcem$ un revêtement universel riemannien de $\srfcem$. Si on munit le bord à l'infini
$\partial_\infty\widetilde{\Sigma}$ de l'ordre $o_\infty$ défini dans la remarque \ref{ordresurface}, on est dans le cadre de la partie \ref{CAT}. On considère alors les
paramétrages (à changements d'origine près) des géodésiques hyperboliques locales de $\Lambda_m$, en considérant pour chacune un paramétrage dans chaque sens. On note 
$\Gamma_\Sigma$ le groupe de revêtement de $p$ et $\widetilde{\Lambda}_{\widetilde{m}}$ l'image réciproque      
    de  $\Lambda_m$ dans $\revetm$. Les géodésiques de $\Lmr$ sont deux à deux disjointes et donc   
       deux à deux non entrelacées. De plus, puisque le support de $\Lmr$ est fermé dans $\widetilde{\Sigma}$ et la métrique est hyperbolique, l'ensemble
des couples de points à l'infini des éléments de $\Lmr$ est fermé dans $\dddp$ (d'après \cite[Lem.~1.6.1~p.68]{Penner92}), et puisque $\revetm$ est $\CAT(-1)$, pour chaque
élément $(x,y)\in\dddp$, il existe une unique géodésique (à changement d'origine près) dont le couple de points à l'infini est $(x,y)$. D'après le lemme \ref{Xfermé}, l'ensemble
$\Lmr$ est fermé pour la topologie des géodésiques. C'est donc une lamination au sens de la définition \ref{laminationplate1}, et puisque il est 
$\Gamma_\Sigma$-équivariant, l'ensemble des projetés des feuilles de $\Lmr$ est une lamination de $\srfcem$.
Dans la suite de cette partie, on appellera $\Lambda_m$ cet ensemble et support de $\Lambda_m$ la réunion $\operatorname{Supp}\Lambda_m$ des images de ses éléments dans $\Sigma$.

%
%

La lamination hyperbolique $\Lambda_m$
est {\it minimale} si elle ne contient
pas de sous-lamination stricte. De manière équivalente, $\Lambda_m$ est 
minimale si l'image de tout rayon géodésique de toute feuille de $\Lambda_m$ est dense dans le support de $\Lambda_m$. 
Une {\it composante minimale} d'une lamination est une sous-lamination minimale.

\medskip

Supposons que $\Sigma$ soit d'aire finie (pour la métrique hyperbolique). Si on munit le complémentaire de $\operatorname{Supp}\Lambda_m$ dans
$\Sigma$ de la distance de longueur induite par la restriction de la métrique hyperbolique sur $\Sigma-\operatorname{Supp}\Lambda_m$, 
le complété $\overline{D}$
d'une composante connexe $D$ du 
complémentaire pour cette distance est une surface hyperbolique à bord totalement géodésique, d'aire finie. Il est obtenu en rajoutant abstraitement à $D$ un bord 
qui est la réunion des images 
d'un nombre fini de feuilles de $\Lambda_m$. En particulier, le bord du complété de $D$ a un nombre fini de pointes, chacune délimitée par deux géodésiques du bord 
(éventuellement égales), voir \cite[p.7]{Bonahon97}.

Si $\srfcem$ est d'aire finie, on dit qu'un bout $e$ d'une géodésique locale simple $\lambda$ de $\srfcem$ {\it spirale} sur une
lamination hyperbolique $\Lambda_m$ si $\lambda$ n'intersecte pas le support de $\Lambda_m$
et s'il existe un rayon $r$ dans la classe de $e$ tel que le support de $\Lambda_m$ soit l'ensemble 
des points d'accumulation de $r(t)$ quand $t$ tend vers $+\infty$. Cette définition ne dépend pas du choix du rayon 
(voir la partie \ref{CAT} pour la définition d'un bout).

\blemm\label{spirale}
Soit $\Lambda_m$ une lamination hyperbolique de $\srfcem$. On considère une feuille $\lambda$ de $\Lambda_m$, un bout $e$ de $\lambda$ 
(correspondant à $t\to+\infty$), et $C$ une 
composante minimale de $\Lambda_m$ qui n'est pas une feuille fermée. Alors $e$ spirale sur $C$ si et seulement s'il existe une feuille $\lambda_0$ de $C$ 
(égale à $\lambda$ si $\lambda$ appartient à $C$)   
dont le point à l'infini (pour $t\to+\infty$) d'un relevé est égal au point à l'infini d'un relevé d'un rayon géodésique local dans la classe de $e$.
\elemm

\dem On fixe un paramétrage de $\lambda$.  Supposons qu'il existe une feuille $\lambda_0$ de $C$ dont un relevé $\widetilde{\lambda}_0$ a 
le même point à l'infini (pour $t\to+\infty$) que celui d'un relevé  $\widetilde{\lambda}$ de $\lambda$. La projection $p:\revetm\to\srfcem$ n'augmente pas 
les distances, et 
par une propriété bien connue de la géométrie hyperbolique, on sait que quitte à changer d'origine 
$d(\widetilde{\lambda}_0(t),\widetilde{\lambda}(t))$ tend vers $0$, et donc $d(\lambda(t),\lambda_0(t))$ tend vers $0$ quand $t$ tend vers l'infini.
Donc l'ensemble des points d'accumulation 
de $\lambda(t)$ quand $t$ tend vers l'infini est égal à celui de $\lambda_0(t)$, qui est le support de $C$ par minimalité.

%

\medskip

Réciproquement, supposons que le support de $C$ soit l'ensemble limite de $\lambda(t)$ (pour $t\to+\infty$). Si $\lambda$ appartient à $C$, on prend $\lambda_0=\lambda$. 
Sinon, l'image de $\lambda$ est contenue dans une composante $D$ du complémentaire du support de
$C$, et on a vu que le complété $\overline{D}$ de $D$, pour la distance induite par $m$, est obtenu en recollant abstraitement un
bord
à $D$ qui est la réunion des images d'un nombre fini de feuilles de $C$. Puisque le bout  $e$ s'accumule sur $C$, la géodésique $\lambda$ 
converge vers une pointe du bord de $\overline{D}$, et donc il existe une feuille $\lambda_0$ de $C$, abstraitement contenue dans le bord de $\overline{D}$ telle que, 
quitte à
changer d'origine, la distance $d(\lambda(t),\lambda_0(t))$ tend vers $0$. Il existe alors des relevés $\widetilde{\lambda}$ et $\widetilde{\lambda}_0$
de $\lambda$ et $\lambda_0$ tels que $d(\widetilde{\lambda}(t),\widetilde{\lambda}_0(t))$ tend vers $0$.\cqfd

\medskip

Soient $\Lm$ une lamination géodésique hyperbolique de $\srfcem$ et $\widetilde{\Lambda}_{\widetilde{m}}$ son image réciproque dans $\widetilde{\Sigma}$.
La lamination $\Lm$ est dite {\it maximale} si ce n'est pas une sous-lamination stricte d'une lamination géodésique hyperbolique.
Si $(\Sigma,m)$ est d'aire finie, la lamination $\Lm$ est maximale si les complétés des composantes connexes du
complémentaire du support de $\Lambda_m$ dans $\Sigma$ (pour la métrique induite) sont des triangles géodésiques idéaux ou des monogones géodésiques 
idéaux privés d'un point.
Si $\srfcem$ est une surface compacte,  à bord (éventuellement vide) totalement géodésique, on dit que
 $\Lm$ est {\it quasi-maximale} si aucune composante de bord n'appartient à $\Lm$ et la réunion de $\Lm$ et des composantes de bord est maximale. La lamination $\Lm$
 est donc quasi-maximale si et seulement si les complétés des composantes connexes du complémentaire du support de la lamination sont des triangles géodésiques 
idéaux ou des monogones géodésiques idéaux privés d'un disque dont le bord n'est pas une feuille. 

\blemm\label{maximaleminimale}
Si $\srfcem$ est une surface hyperbolique, connexe, compacte, à bord totalement géodésique, il existe au moins une lamination sur $\srfcem$ qui est à la fois minimale et quasi-maximale si et seulement si $\Sigma$ n'est pas un pantalon.
\elemm

\dem Soient $(\Sigma_1,m_1)$ la surface hyperbolique complète obtenue en recollant des bouts évasés à bord géodésique sur chacune des composantes de bord de $\Sigma$
et $m_1'$ une métrique hyperbolique sur $\Sigma_1$, complète et d'aire finie. Puisque $\Sigma$
n'est pas un pantalon, la surface $\Sigma_1$  n'est pas la sphère privée de trois points. Il existe alors une lamination $\Lambda_1$ de $(\Sigma_1,m_1')$
qui est à la fois minimale et maximale, dont les feuilles ne finissent pas dans les pointes (voir par exemple la démonstration de \cite[Lem.~2.3~p.~10-11]{Hamenstadt09a}).
En étendant naturellement la bijection naturelle de
\cite[Lem.~18]{Bonahon97} entre laminations géodésiques pour des métriques hyperboliques qui ne sont pas nécessairement d'aire finie, et dont les feuilles ne finissent
pas dans les pointes ou les bouts évasés, on montre qu'il existe une unique lamination géodésique de $\srfcem$ en bijection avec $\Lambda_1$, qui est minimale comme
$\Lambda_1$ et qui est quasi-maximale car $\Lambda_1$ est maximale.  Enfin, si $\Sigma$ est un pantalon, les seules laminations minimales sont  les géodésiques du 
bord, qui ne sont pas quasi-maximales.\cqfd

\subsection{Liens entre géodésiques locales hyperboliques et plates.}\label{hyperboliqueplate}

On suppose dans cette partie \ref{hyperboliqueplate} que $\Sigma$ est compacte et $\chi(\Sigma)<0$. On considère une structure de demi-translation $[q]$ ainsi
qu'une métrique hyperbolique $m$ (complète, à bord totalement géodésique) sur $\Sigma$, et on fixe 
$p:\widetilde{\Sigma}
\to\Sigma$ un revêtement universel de $\Sigma$. On note $\Gamma_\Sigma$ le groupe de revêtement de $p$, $[\widetilde{q}]$ l'unique structure de demi-translation et $\widetilde{m}$ l'unique métrique hyperbolique 
définies
sur $\widetilde{\Sigma}$ telles que $p:(\widetilde{\Sigma},[\widetilde{q}])\to(\Sigma,[q])$ et $p:(\widetilde{\Sigma},\widetilde{m})\to(\Sigma,m)$ soient
localement isométriques. Les espaces métriques
$\revet$ et $\revetm$ sont $\CAT(0)$ et puisque $\Sigma$ est compacte, il existe un 
unique homéomorphisme $\Ga_\Sigma$-équivariant entre les bords à l'infini de $\widetilde{\Sigma}$ pour les métriques $\widetilde{m}$ et $[\widetilde{q}]$,
que l'on note indistinctement
$\partial_\infty\widetilde{\Sigma}$, et $\partial_\infty\Gamma_\Sigma$, par lequel on les identifie (voir par exemple \cite[§2]{Bonahon91}). 

Si $G$ est un ensemble de géodésiques hyperboliques de $\widetilde{\Sigma}$ deux à deux non entrelacées, alors $\operatorname{Supp} G$ est fermé dans $\widetilde{\Sigma}$ 
si et seulement si $E(G)$ est fermé dans $\dddp$ (voir \cite[Lem.~1.6.1~p.~68]{Penner92} et la partie \ref{CAT} pour les définitions de $\operatorname{Supp} G$ et $E(G)$).
Dans le cas des ensembles de géodésiques plates, cette équivalence est fausse.
Dans cette partie, nous étudions le lien entre ensembles de géodésiques plates et ensembles
de géodésiques hyperboliques.

\medskip

On appelle {\it lamination plate} une lamination géodésique au sens de la définition \ref{laminationplate1} d'une surface munie d'une structure de demi-translation.
On dit qu'elle est {\it pleine} si son image réciproque (par $p$) n'est pas une sous-lamination stricte d'une lamination plate ayant le même ensemble de couples
de points à l'infini de feuilles. Toute lamination plate $\Lq$ est contenue dans une unique lamination plate pleine $\Lq^{\pl}$ ayant le même ensemble de couples
de points à 
l'infini au revêtement universel et $\Lq^{\pl}-\Lq$ est composée de feuilles périodiques (voir le lemme \ref{geodasymptotes}).
On dit que $\Lq^{\pl}$ est le {\it remplissage} de $\Lq$.  
Si $\lambda$ et $\ell$ sont des géodésiques locales (non paramétrées mais orientées) de $\srfcem$ et $\srfce$, on dit que $\lambda$ {\it correspond} à 
$\ell$ s'il existe des relevés $\widetilde{\lambda}$ et $\widetilde{\ell}$ de 
$\lambda$ et $\ell$ ayant le même couple de points à l'infini, et on dit qu'une lamination hyperbolique {\it correspond} à une lamination plate
si leurs images réciproques au revêtement universel ont les mêmes ensembles de couples de points à l'infini de feuilles. À toute géodésique locale plate 
correspond une unique géodésique locale hyperbolique. De plus, si $\widetilde{\Lambda}_{[\widetilde{q}]}$ est l'image réciproque par $p$ d'une lamination plate
$\Lq$, et si $Y$ est l'ensemble des couples de points à l'infini des feuilles de $\Lqr$, alors $Y$ est 
$\Gamma_\Sigma$ et $\iota$-invariant (où $\iota:(x,y)\mapsto (y,x)$), les couples d'éléments de $Y$ sont deux à deux non entrelacés et  $Y$ est fermé dans 
$\dddp$  d'après le lemme
\ref{caractérisation}. D'après \cite[Lem.~1.6.1~p.68]{Penner92}, il existe donc une unique lamination hyperbolique $\Lmr$ de 
$\revetm$ telle que 
$E(\Lmr)=Y$, et $\Lmr$ est $\Gamma_\Sigma$-invariante par naturalité. Donc l'ensemble
 $\Lm$ des projetés des feuilles de $\Lmr$ est l'unique lamination hyperbolique de $\srfcem$ qui correspond à $\Lq$ (par construction).

\blemm\label{bijectionplatmaxhyperbolique}
L'application qui à une lamination plate associe l'unique lamination hyperbolique qui lui correspond induit une bijection entre l'ensemble des laminations
 plates pleines et l'ensemble des laminations hyperboliques.
\elemm

\dem Montrons que la restriction de cette application aux laminations plates pleines est surjective. Soient $\Lambda_m$ une lamination hyperbolique de
$\srfcem$, et $\widetilde{\Lambda}_{\widetilde{m}}$ son image réciproque par $p$.
On note $Y=E(\Lmr)$,   
$\widetilde{\Lambda}^{\pl}_{[\widetilde{q}]}$ l'ensemble maximal de géodésiques de $\revet$ (définies à changements d'origine près) tel que pour tout 
$\widetilde{\ell}\in\widetilde{\Lambda}^{\pl}_{[\widetilde{q}]}$ il existe $\widetilde{\lambda}\in\widetilde{\Lambda}_{\widetilde{m}}$ qui correspond à 
$\widetilde{\ell}$.     

%

 Comme pour $\widetilde{\Lambda}_{\widetilde{m}}$, les géodésiques de 
$\widetilde{\Lambda}^{\pl}_{[\widetilde{q}]}$ ne sont pas entrelacées et $\widetilde{\Lambda}^{\pl}_{[\widetilde{q}]}$ est stable par $\ell\mapsto\ell^- $.
         Enfin, puisque $Y$ est fermé dans $\dddp$, et puisque l'ensemble $\widetilde{\Lambda}^{\pl}_{[\widetilde{q}]}$ contient toutes les géodésiques plates dont le couple de points
à l'infini appartient à $Y$, d'après le lemme \ref{Xfermé},   
 l'ensemble $\widetilde{\Lambda}^{\pl}_{[\widetilde{q}]}$ est fermé pour la topologie des géodésiques, et c'est donc une lamination plate pleine. 
 Et comme $\widetilde{\Lambda}_{\widetilde{m}}$, l'ensemble $\widetilde{\Lambda}^{\pl}_{[\widetilde{q}]}$ est 
$\Gamma_\Sigma$-équivariant. Donc l'ensemble $\Lambda_{[q]}^{\pl}$ des projetés des feuilles de $\widetilde{\Lambda}^{\pl}_{[\widetilde{q}]}$ par $p$ 
est une lamination plate de $\srfce$ et, par construction,  $\Lm$ correspond à $\Lambda^{\pl}_{[q]}$.
Donc la restriction de l'application aux laminations plates pleines est surjective. Par définition des laminations plates pleines, elle est aussi
injective, et donc bijective.\cqfd
%

%

%
%
%

\medskip
On rappelle le théorème de la bande plate. On dit que deux géodésiques $c$ et $c'$ d'un espace métrique $(X,d)$ sont \textit{à distance de Hausdorff finie}
s'il existe $K>0$ tel que
$d(c(t),c'(t))\leqslant K$ pour tout $t\in\RR$.

\btheo(voir \cite[Th.~2.13~p.182]{BriHae99})\label{bandeplate}
Soient $(X,d)$ un espace métrique $\CAT(0)$ et $c,c':\RR\to X$ deux géodésiques. Si $c$ et $c'$ sont à distance de Hausdorff finie, l'enveloppe convexe de 
$c(\RR)\cup c'(\RR)$ est isométrique
à une bande plate $\RR\times[0,D]\subset \EE^2$, avec $D\geqslant 0$.\cqfd 
\etheo

%

\blemm(voir \cite[Th.~2.(c)]{MS85})\label{geodasymptotes}
On rappelle que $\Sigma$ est compacte. Soient $\widetilde{\ell}_1$ et $\widetilde{\ell}_2$ deux géodésiques de $\revet$, telles que $\widetilde{\ell}_1$ et $\widetilde{\ell}_2$ sont à distance de Hausdorff
finie et qu'il existe $\delta>0$ tel que 
$d(\widetilde{\ell}_1(t),\widetilde{\ell}_2(\RR))\geqslant\delta$ pour tout $t\in\RR$. Alors les projetés $p\circ\widetilde{\ell}_1$ et $p\circ\widetilde{\ell}_2$
de $\widetilde{\ell}_1$
et $\widetilde{\ell}_2$ sur $\Sigma$ sont des géodésiques périodiques qui sont librement homotopes au bord d'un même cylindre plat maximal, dont l'intérieur est
plongé isométriquement dans $\srfce$.\cqfd
\elemm

\bcoro\label{corogeodasymptotes}
Soit $\Lambda_m$ une lamination hyperbolique de $(\Sigma,m)$ qui ne contient pas de feuille fermée. Alors il existe une unique lamination plate 
$\Lambda_{[q]}$ telle que $\Lm$ corresponde à $\Lq$. Elle est pleine et à chaque feuille de $\Lambda_m$ ne correspond qu'une seule feuille de $\Lambda_{[q]}$.
\ecoro

\dem Supposons qu'il existe deux géodésiques (non paramétrées) distinctes $\widetilde{\ell}_1$ et $\widetilde{\ell}_2$ de $\revet$ qui ont le même couple de points à
l'infini qu'un relevé $\widetilde{\lambda}$ d'une feuille $\lambda$ de $\Lambda_m$. Alors d'après le théorème \ref{bandeplate}, les géodésiques
$\widetilde{\ell}_1$ et
$\widetilde{\ell}_2$ bordent une même bande plate, et il existe $\delta>0$ tel que tout point de l'une est à distance $\delta$ de l'autre. Donc, d'après le lemme 
\ref{geodasymptotes}, les projetés $\ell_1=p\circ\widetilde{\ell}_1$ et $\ell_2=p\circ\widetilde{\ell}_1$ sont 
périodiques.  Mais alors il existe un élément hyperbolique $\gamma$ du groupe de revêtement de $p$, tel que 
$\widetilde{\lambda}(-\infty)=\widetilde{\ell}_1(-\infty)$ et
$\widetilde{\lambda}(+\infty)=\widetilde{\ell}_1(+\infty)$ sont respectivement les points fixes répulsif et attractif de $\gamma$, et 
 $\lambda$ serait fermée. Donc, à toute feuille de $\Lambda_m$ ne correspond qu'une
seule géodésique locale de $\srfce$, et il n'existe donc qu'une seule lamination plate qui correspond à $\Lambda_m$. 
\cqfd



\subsection{Laminations plates minimales.}\label{partiecomposanteminimale}

Dans cette partie \ref{partiecomposanteminimale}, on suppose toujours que $\Sigma$ est compacte et $\chi(\Sigma)<0$ et on note $p:\revet\to\srfce$ un revêtement localement
isométrique de groupe de revêtement $\Gamma_\Sigma$. De même que pour une lamination hyperbolique,
une lamination plate est dite {\it minimale} si elle n'a pas de sous-lamination stricte. De manière équivalente, une lamination plate est minimale si pour
chacune de ses feuilles $\ell$, la paire $\{\ell,\ell^-\}$, où $\ell^-$  est la feuille $\ell$ orientée dans l'autre sens, est dense dans la lamination, 
pour la topologie des géodésiques. Si une lamination plate de $\srfce$ est minimale,
la lamination hyperbolique de $\srfcem$ qui lui correspond n'a pas de sous-lamination stricte et est donc minimale. Réciproquement, supposons 
qu'une lamination hyperbolique est minimale et n'est
pas une feuille fermée, alors par minimalité, elle ne contient pas de feuille fermée et d'après le corollaire \ref{corogeodasymptotes}, 
il existe une unique lamination plate qui lui correspond, qui est minimale.  On considère deux lemmes généraux.

\blemm\label{minimalepaspériodique}
Si $\Lambda$ est une lamination plate minimale de $\srfce$, qui n'est pas une paire de feuilles périodiques opposées,
aucune de ses feuilles n'est positivement ou négativement périodique.
\elemm

\dem Supposons par contraposée qu'il existe une feuille $\ell$ de $\Lambda$ qui est positivement périodique. Alors quitte à changer son origine, elle est égale, à partir
d'un certain temps, à une géodésique locale périodique $\ell'$ et il existe une suite de réels $(t_n)_{n\in\NN}$ tendant vers l'infini telle que 
$\ell(t_n)=\ell'(0)$ pour tout $n$. Alors la suite $(\ell_n)_{n\in\NN}$
définie par $\ell_n(t)=\ell(t+t_n)$ pour tout $t\in\RR$ converge vers $\ell'$ pour la topologie compacte-ouverte, et $\ell'$ appartient à $\Lambda$ car $\Lambda$ est fermée pour 
la topologie des géodésiques. Par minimalité, $\Lambda=\{\ell',\ell'^-\}$.\cqfd

\medskip

Soient $\Lambda$ une lamination hyperbolique de $\srfcem$ et $C$ une composante minimale de $\Lambda$ qui n'est pas une feuille fermée.
Si $\lambda_2$ est
une feuille de $C$ et si $\lambda_1$ est une feuille de $\Lambda$ dont le bout (pour $t\to+\infty$) spirale sur $C$, il existe des géodésiques
locales plates
$\ell_1$ et $\ell_2$, uniques à changements d'origine près, qui correspondent à $\lambda_1$ et $\lambda_2$. On fixe des paramétrages de $\lambda_1$, $\lambda_2$,
$\ell_1$ et $\ell_2$. 

\blemm\label{composanteminimale}
Pour tout relevé paramétré $\widetilde{\ell}_2$ de $\ell_2$ dans $\widetilde{\Sigma}$, quitte à inverser l'orientation de $\ell_2$,
il existe une suite de réels $(t_n)_{n\in\NN}$ qui converge vers l'infini, et une suite de relevés (paramétrés) $(\widetilde{\ell}_{1,n})$ de $\ell_1$ tels que
la suite de feuilles paramétrées $(\widetilde{\ell}_{1,n}')_{n\in\NN}$ définies par $\widetilde{\ell}_{1,n}'(t)=\widetilde{\ell}_{1,n}(t+t_n)$ pour tout $t\in\RR$, converge vers
 $\widetilde{\ell}_2$ pour la topologie compacte-ouverte. 
 \elemm


\dem Soit $\widetilde{\lambda}_2$ le relevé de $\lambda_2$ qui correspond à $\widetilde{\ell}_2$.
Puisque le bout de $\lambda_1$ (pour $t\to+\infty$)
spirale sur $C$, 
il existe une suite de réels $(s_n)_{n\in\NN}$ tendant vers l'infini telle que $\lambda_1(s_n)$ converge vers $\lambda_2(0)$, et il existe une suite
$(\widetilde{\lambda}_{1,n})_{n\in\NN}$ de relevés de $\lambda_1$ telle que la suite $(\widetilde{\lambda}_{1,n}(s_n))_{n\in\NN}$ converge vers
$\widetilde{\lambda}_2(0)$. D'après le lemme \ref{compac}, si pour tout $n\in\NN$, on reparamètre $\widetilde{\lambda}_{1,n}$ pour que 
$\widetilde{\lambda}_{1,n}(s_n)$ soit l'origine, quitte à extraire et à inverser l'orientation de $\lambda_2$,
la suite $(\widetilde{\lambda}_{1,n})_{n\in\NN}$ converge vers $\widetilde{\lambda}_2$ pour la topologie compacte-ouverte, et donc la suite des couples de points à l'infini des 
feuilles $(\widetilde{\lambda}_{1,n})_{n\in\NN}$ converge vers le couple de points à l'infini de $\widetilde{\lambda}_2$. 
Soient $(\widetilde{\ell}_{1,n})_{n\in\NN}$ la suite de relevés de $\ell_1$ telle que, pour tout $n\in\NN$, $\widetilde{\ell}_{1,n}$
corresponde à $\widetilde{\lambda}_{1,n}$ et $(t_n)_{n\in\NN}$ la suite de réels telle que $\widetilde{\ell}_{1,n}(t_n)$ soit un 
point de $\widetilde{\ell}_{1,n}(\RR)$
le plus proche de $\widetilde{\ell}_{2}(0)$. Puisque $\Sigma$ est compacte, il existe $\delta>0$ tel que $\revet$ est $\delta$-hyperbolique, 
et d'après le lemme \ref{convergenceinfini} et la définition de $(t_n)_{n\in\NN}$, quitte à extraire,
la suite $(\widetilde{\ell}_{1,n}')_{n\in\NN}$ définie par 
$\widetilde{\ell}_{1,n}'(t)=\widetilde{\ell}_{1,n}(t+t_n)$ pour tout $t\in\RR$, converge vers
$\widetilde{\ell}_2$ pour la topologie compacte-ouverte.

De plus, pour tout $n\in\NN$, il existe $\gamma_n\in\Gamma_\Sigma$ (le groupe de revêtement de $p$) tel que
$\widetilde{\ell}_{1,n}=\gamma_n\widetilde{\ell}_{1,1}$, et puisque $\widetilde{\ell}_{1,n}$ et $\widetilde{\lambda}_{1,n}$ ont le même couple de points à l'infini,
 on a $\widetilde{\lambda}_{1,n}=\gamma_n\widetilde{\lambda}_{1,1}$. Puisque $\Gamma_\Sigma$ agit par isométrie sur $\revet$, on en déduit que 
 $d(\widetilde{\ell}_{1,n}(0),\widetilde{\lambda}_{1,n}(0))=d(\widetilde{\ell}_{1,1}(0),\widetilde{\lambda}_{1,1}(0))$, pour la distance définie par
 la structure de demi-translation.
 De même, les points $\widetilde{\ell}_{1,n}(t_n)$ et
 $\widetilde{\lambda}_{1,n}(s_n)$ restent à distance bornée. De plus, puisque $\Sigma$ est compacte, les revêtements $\revet$ et $\revetm$ sont quasi-isométriques, donc les
 géodésiques hyperboliques $\widetilde{\lambda}_{1,n}$ sont des quasi-géodésiques de $\revet$ (voir par exemple \cite[Prop.~8.19 et Déf.~8.22 p.~140-142]{BriHae99}).  
 Donc la distance $d(\widetilde{\lambda}_{1,n}(s_n),\widetilde{\lambda}_{1,n}(0)$ tend vers l'infini et donc la suite $(t_n)_{n\in\NN}$ tend vers l'infini.\cqfd


\bcoro\label{corocomposanteminimale}
Soit $\Lq$ une lamination plate minimale qui n'est pas une paire de feuilles périodiques opposées. Alors 
si $\ell_1$ et $\ell_2$ sont deux feuilles de $\Lambda_{[q]}$ et $\widetilde{\ell}_2$ est un relevé de $\ell_2$ dans $\widetilde{\Sigma}$, 
quitte à inverser l'orientation
de $\widetilde{\ell}_2$, il existe une suite de relevés de 
$\ell_1$ qui converge vers $\widetilde{\ell_2}$ pour la topologie des géodésiques. De plus, l'image de chacun des rayons géodésiques de chacune des feuilles de 
$\Lq$ est dense dans le support de $\Lq$.
\ecoro

\dem On applique le lemme \ref{composanteminimale} en prenant $C=\Lambda_m$ qui est la lamination hyperbolique correspondant à $\Lq$, qui est minimale et n'est pas une
feuille fermée, en remarquant que chacun des bouts de chacune
des feuilles de $\Lm$ spirale sur $\Lm$. Soient $\ell_1$ et $\ell_2$ deux feuilles de $\Lq$ et $x$ un point de l'image de $\ell_2$. On fixe des paramétrages
tels que $\ell_2(0)=x$. Alors il existe une suite $(t_n)_{n\in\NN}$ qui tend vers l'infini et une suite $(\widetilde{\ell}_{1,n})_{n\in\NN}$ de relevés de $\ell_1$
telles que $(\widetilde{\ell}_{1,n}(t_n))_{n\in\NN}$ converge vers $\widetilde{\ell}_2(0)$ et, quitte à changer l'orientation de $\ell_2$ et à prendre
$\widetilde{\ell}_{1,n}(t_n)$ comme origine de $\widetilde{\ell}_{1,n}$ pour tout $n\in\NN$, la suite 
$(\widetilde{\ell}_{1,n})_{n\in\NN}$ converge vers $\widetilde{\ell_2}$ pour la topologie compacte-ouverte. Donc pour tout $\varepsilon>0$ et $t>0$, il existe 
$n\in\NN$ tel que $t_n>t$ et $d(\widetilde{\ell}_{1,n}(t_n),\widetilde{\ell}_2(0))=d({\ell}_{1}(t_n),x)<\varepsilon$  
(pour la distance de la structure de demi-translation). Donc les images de tous rayons de la forme $\ell_{1|[a,+\infty[}$ avec $a\in\RR$ sont denses dans le support de $\Lq$ et de mêmes pour 
les rayons de la forme $\ell_{1|]-\infty,a]}$. \cqfd

\begin{center}
  \begin{picture}(0,0)%
\includegraphics{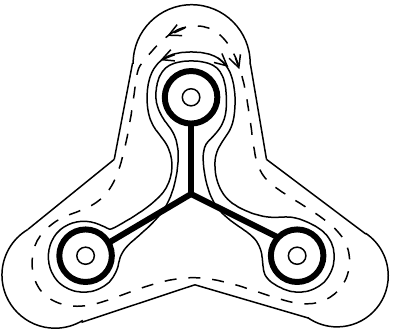}%
\end{picture}%
\setlength{\unitlength}{3771sp}%
\begingroup\makeatletter\ifx\SetFigFont\undefined%
\gdef\SetFigFont#1#2#3#4#5{%
  \reset@font\fontsize{#1}{#2pt}%
  \fontfamily{#3}\fontseries{#4}\fontshape{#5}%
  \selectfont}%
\fi\endgroup%
\begin{picture}(1959,1644)(4755,-2949)
\end{picture}%

 \end{center}
 
\rem Une lamination plate peut ne pas être minimale même si les images de tous les rayons géodésiques de toutes ses feuilles sont denses dans
son support. Par exemple, on considère la surface de demi-translation ci-dessus, dont les singularités sont d'angle $3\pi$ et situées aux sommets
du graphe en gras. Alors la réunion des deux géodésiques plates correspondant aux deux courbes périodiques dessinées (et de leurs inverses) 
est une lamination plate telle que  l'image de chacun des rayons géodésiques de chacunes des feuilles est égale au graphe. Mais la lamination n'est pas minimale.

Dans la suite de cette partie \ref{partiecomposanteminimale}, on considère une lamination plate minimale $\Lq$ de $\srfce$, qui n'est pas une paire de feuilles
périodiques opposées et 
la lamination hyperbolique minimale $\Lm$ de $\srfcem$ qui lui correspond, qui n'est pas une feuille fermée. On note 
$\widetilde{\Lambda}_{[\widetilde{q}]}$ et $\widetilde{\Lambda}_{\widetilde{m}}$ leurs images réciproques dans $\widetilde{\Sigma}$.

\blemm\label{composanteminimalecompacte}
 Soit $\Lq$ une lamination plate minimale de $\srfce$ qui n'est pas une paire de feuilles périodiques opposées.  Supposons qu'il existe une feuille de $\Lambda_{[q]}$ dont l'image est compacte. Alors toutes les feuilles de $\Lambda_{[q]}$ ont la même image, 
et l'image de n'importe quel rayon géodésique contenu dans une feuille lui est égale.  
\elemm

\dem Soient $\ell_1,\ell_2$ des feuilles de $\Lambda_{[q]}$. D'après le corollaire \ref{corocomposanteminimale},
si $\widetilde{\ell}_2$
 est un relevé de $\ell_2$, il existe une suite de relevés de $\ell_1$ qui converge vers $\widetilde{\ell_2}$ pour la topologie des géodésiques.

Supposons que l'image de $\ell_1$ est compacte. Si $\ell_1$ était régulière, son image serait contenue dans un cylindre plat et $\ell_1$ serait périodique, et on aurait $\Lq=\{\ell_1,\ell_1^-\}$. 
Donc
l'image de $\ell_1$ est contenue dans une union finie de liaisons de singularités. Si $K$  est un compact
de $\widetilde{\Sigma}$ qui intersecte l'image de $\widetilde{\ell}_2$, il n'y a qu'un nombre fini de liaisons de singularités contenues dans $K$ qui se projettent sur
une liaison de singularités parcourue par $\ell_1$. Et puisque $\revet$ est un espace métrique $\CAT(0)$, tout segment géodésique parcourt au plus une fois chacune de
ces liaisons de singularités. Il n'y a donc qu'un nombre fini de segments géodésiques d'image contenue dans $K$ dont l'image est contenue dans la réunion
de ces liaisons de singularités, et la suite $(\widetilde{\ell}_{1,n|K})_{n\in\NN}$ des restrictions des feuilles $\widetilde{\ell}_{1,n}$ à $K$ est constante
à partir d'un certain rang, égale à $\widetilde{\ell}_{2|K}$. 
On en déduit que pour tout compact 
$K$ de $\widetilde{\Sigma}$, l'image du projeté $p(\widetilde{\ell}_{2|K})$ est contenue dans l'image de $\ell_1$, et ceci étant vrai pour 
tous les compacts, l'image de $\ell_2$ est contenue dans celle de $\ell_1$, et en particulier compacte. En inversant les rôles de $\ell_1$ et $\ell_2$, on montre que 
l'image de $\ell_1$ est contenue dans celle de $\ell_2$, et donc, elles sont égales. Toutes les feuilles de $\Lambda_{[q]}$ ont donc la même image qui est 
la réunion d'un nombre fini de liaisons de singularités. De plus, d'après le lemme \ref{composanteminimale}, chacun des rayons géodésiques contenus dans une
feuille de $\Lq$ est dense dans le support de $\Lq$, et donc égal au support de $\Lq$.\cqfd

\medskip

Si $\widetilde{\ell}$ est une géodésique de $\revet$, on dit que $\widetilde{\ell}$ est une {\it feuille de Levitt} s'il existe
une composante connexe $\widetilde{\Sigma}^+$ de 
$\widetilde{\Sigma}-\widetilde{\ell}(\RR)$ 
telle que pour tout 
$t\in\RR$, l'angle
défini par les germes de $\widetilde{\ell}(]t-\varepsilon,t])$ et
$\widetilde{\ell}([t,t+\varepsilon[)$, avec $\varepsilon>0$, mesuré dans $\widetilde{\Sigma}^+$, est égal à $\pi$
 (la notion de feuille de Levitt non régulière correspond aux droites singulières du 
relevé d'un feuilletage introduites par G.~Levitt dans \cite{Levitt81}). On appelle $\widetilde{\Sigma}^+$ un {\it côté sans angle} de $\widetilde{\ell}$. 
 On appelle {\it feuille de Levitt} de $\srfce$ le projeté d'une feuille
de Levitt de $\revet$. Une feuille de Levitt est de direction constante.

\medskip

Supposons qu'il existe une feuille $\ell$ de $\Lambda_{[q]}$  dont l'image n'est pas compacte. Puisque $\ell$ n'est pas périodique, d'après les lemmes
\ref{directionconstante}
et \ref{structuredestrajectoires}, elle admet un rayon géodésique $\ell_{|[T,+\infty[}$, avec $T>0$, ne rencontrant pas de singularité, qui est de direction
constante et est dense dans un domaine du feuilletage
vertical $\F_{q'}$ d'une différentielle quadratique holomorphe $q'$ de la classe $[q]$, que l'on note $D$. D'après le corollaire \ref{corocomposanteminimale}, 
tous les rayons géodésiques des feuilles de $\Lq$ sont denses dans le support de $\Lq$, qui est donc égal à $D$. Puisque $\ell([T,+\infty[)$ est dense dans $D$ et
 les autres feuilles de $\Lq$ n'intersectent pas transversalement $\ell(]T,+\infty[)$, toutes les feuilles de $\Lq$ sont de direction constante et de même direction.

\blemm\label{composanteminimalerécurrente}
On rappelle que $\Lq$ est une lamination plate minimale qui n'est pas une paire de feuilles
périodiques opposées, et contient une feuille $\ell$ qui n'est pas d'image compacte. La lamination plate $\Lambda_{[q]}$ est alors égale à
l'ensemble des feuilles de Levitt qui 
sont contenues dans $D$ et de même direction que $\ell$.
\elemm


\dem Le support de la lamination plate $\Lambda_{[q]}$ est $D$ donc $\Lq$ contient l'ensemble des feuilles régulières du feuilletage $\F_{q'}$ 
(le feuilletage vertical défini par $q'$) de même direction que $\ell$, qui sont denses dans $D$.


Montrons que l'ensemble des feuilles de Levitt locales, contenues dans $D$ et de même direction que $\ell$, est contenu dans $\Lambda_{[q]}$. 
Soient $\ell_1$ une feuille de Levitt contenue dans
$D$ et de même direction que $\ell$ et $\widetilde{\ell}_1$ un de ses relevés.     

Soient $\widetilde{\Sigma}^+$ un des côtés sans angle de $\widetilde{\ell}_1$, $\widetilde{x}\in\widetilde{\ell}_1(\RR)$
et $(\widetilde{x}_n)_{n\in\NN}$
une suite de points telle que pour tout $n\in\NN$, il existe une feuille régulière $\widetilde{\ell}_n$ de $\Lqr$ 
contenues dans $\widetilde{\Sigma}^+$ telle que
$\widetilde{\ell}_n(0)=\widetilde{x}_n$, et 
$(\widetilde{x}_n)_{n\in\NN}$ converge vers $\widetilde{x}$.
D'après le lemme \ref{compac}, quitte à extraire, la suite $(\widetilde{\ell}_n)_{n\in\NN}$ converge vers une feuille $\widetilde{\ell}_2$ de $\Lqr$ passant par $\widetilde{x}$,
pour la topologie des géodésiques. Or, pour tout $n\in\NN$, $\widetilde{\ell}_n$ est régulière et contenue dans 
$\widetilde{\Sigma}^+$, donc la géodésique $\widetilde{\ell}_2$ est une feuille de Levitt dont $\widetilde{\Sigma}^+$ est un côté sans angle, et donc 
$\widetilde{\ell}_2=\widetilde{\ell}_1$
à changement d'orientation et d'origine près. Donc $\widetilde{\ell}_1$ appartient à $\widetilde{\Lambda}_{[\widetilde{q}]}$. 

\medskip

Réciproquement, montrons que $\Lambda_{[q]}$ est contenu dans l'ensemble des feuilles de Levitt contenues dans $D$ et de même direction que $\ell$. Nous avons
vu que l'ensemble des feuilles de Levitt régulières, de même direction que $\ell$, qui sont denses dans $D$, est contenu dans $\Lambda_{[q]}$. 
Soient $\ell_1$ l'une d'elles et $\ell_2$
une feuille quelconque de $\Lambda_{[q]}$. Si $\widetilde{\ell}_2$ est un relevé de $\ell_2$, puisque $\Lq$ est minimale et n'est pas réduite à
une paire de feuilles 
périodiques opposées, d'après le corollaire \ref{corocomposanteminimale}, il existe une suite de relevés paramétrés $(\widetilde{\ell}_{1,n})_{n\in\NN}$ de $\ell_1$ qui
converge
vers $\widetilde{\ell}_2$ pour la topologie des géodésiques.
Puisque pour tout $n\in\NN$, les feuilles $\widetilde{\ell}_{1,n}$ et $\widetilde{\ell}_2$ ne sont pas entrelacées, quitte à extraire, il existe une 
composante connexe 
$\widetilde{\Sigma}^+$ de 
$\widetilde{\Sigma}-\widetilde{\ell}_2(\RR)$ dont l'adhérence contient toutes les feuilles de la suite.
Puisque $(\widetilde{\ell}_{1,n})_{n\in\NN}$ converge vers
$\widetilde{\ell}_2$ pour la topologie des géodésiques, pour tout $t\in\RR$, l'angle défini par les germes $\widetilde{\ell}_2(]t-\varepsilon,t])$ et 
$\widetilde{\ell}_2([t,t+\varepsilon[)$, avec $\varepsilon>0$,  
mesuré dans $\widetilde{\Sigma}^+$, est égal à $\pi$. Donc $\widetilde{\ell}_2$ est une feuille de Levitt contenue dans $p^{-1}(D)$ et 
$\ell_2$ est une feuille de Levitt contenue dans $D$.\cqfd
 
\subsection{Rayons asymptotiques.}\label{asymptotes}
%
%

Dans cette partie \ref{asymptotes}, on suppose toujours que $\srfce$ est une surface compacte, à bord (éventuellement vide), munie d'une structure de demi-translation 
et $p:\revet\to\srfce$ est un revêtement universel localement isométrique. On caractérise maintenant les rayons géodésiques de $\revet$ ayant le même point à l'infini.
La notion de rayon géodésique quasi-droit vient de $\cite[\S~3.4]{Dank10}$.

  Soit $\widetilde{r}:\RR^+\to\widetilde{\Sigma}$ un rayon géodésique de $\revet$. On peut prolonger (de manière pas nécessairement unique et indifférente)
  $\widetilde{r}$ en une géodésique. 
On note $\widetilde{\Sigma}^\pm$ les adhérences des composantes connexes du complémentaire de l'image  de cette géodésique dans $\widetilde{\Sigma}$.  
  Pour tout $t>0$,
on note $\theta_{\widetilde{r}}^-(t)$ et $\theta_{\widetilde{r}}^+(t)$ (ou simplement $\theta^-(t)$ et $\theta^+(t)$)   
     les angles définis par les germes de ${\widetilde{r}(]t-\varepsilon,t])}$ et $\widetilde{r}([t,t+\varepsilon[)$, avec $\varepsilon>0$,
     mesurés respectivement dans $\widetilde{\Sigma}^-$ et
$\widetilde{\Sigma}^+$.
On dit que $\widetilde{r}$ est \textit{quasi-droit} si, quitte à permuter $+$ et $-$, on a
$\sum_{t>0}(\theta^+(t)-\pi)<+\infty$.

 \begin{center}
  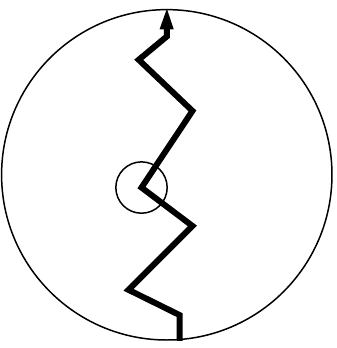
 \end{center}

     Puisque $\widetilde{r}$ est  un rayon géodésique, on a toujours $\theta^+(t)-\pi\geqslant 0$ et $\theta^+(t)-\pi=0$
en dehors des singularités. Un point $\eta\in\partial_\infty\widetilde{\Sigma}$ est
\textit{quasi-droit} s'il existe un rayon géodésique quasi-droit $\widetilde{r}$ telle que $\widetilde{r}(+\infty)=\eta$, et une géodésique est {\it quasi-droite} si tous ses
rayons sont quasi-droits.


\blemm (voir \cite[Prop.~3.10]{Dank10})\label{Dank10}

\noindent
$\bullet$~ Soient $\widetilde{r}_1$ et $\widetilde{r}_2$ deux rayons géodésiques de $\revet$ telles que $\widetilde{r}_1(+\infty)=\widetilde{r}_2(+\infty)$. 
 Alors si l'un n'est pas quasi-droit, quitte à inverser $\widetilde{r}_1$ et $\widetilde{r}_2$,  il existe $T\geqslant 0$ tel que  $\widetilde{r}_1(t)=\widetilde{r}_2(t+T)$ pour tout $t\geqslant 0$.
 
\noindent
$\bullet$~Soit $\eta\in\partial_\infty\widetilde{\Sigma}$ un point quasi-droit. Alors tout rayon géodésique $\widetilde{r}$ de $\revet$ tel
que $\widetilde{r}(+\infty)=\eta$ est quasi-droit. 
\elemm




\blemm\label{spiralesurcompact}
Soient $\ell_1:\RR\to\Sigma$ une géodésique locale positivement compacte de $\srfce$ dont l'image n'est pas contenue dans un cylindre plat non dégénéré
et $\widetilde{\ell}_1$
un relevé de
$\ell_1$ dans $\widetilde{\Sigma}$. Soit $\widetilde{\ell}_2$ une géodésique de $\revet$ 
telle que $\widetilde{\ell}_2(+\infty)=\widetilde{\ell}_1(+\infty)$. Alors les géodésiques
$\widetilde{\ell}_{2}$ et $\widetilde{\ell}_{1}$ coïncident à partir d'un certains temps, quitte à  
changer les origines, et donc $\ell_{2}=p\circ\widetilde{\ell}_{2}$ coïncide aussi avec $\ell_{1}$ à partir de ce temps. 
\elemm

\dem Puisque $\ell_1$ est positivement compacte,  il existe $T_0>0$ tel que les ensembles 
$\{\theta_{\widetilde{\ell}_1}^+(t),t>T_0\}$ et $\{\theta_{\widetilde{\ell}_1}^-(t),t>T_0\}$, définis pour un choix de convention $\pm$,   
       ne prennent qu'un nombre fini de
valeurs qui apparaissent chacune une infinité de fois dans les sommes
$\sum_{t>T_0}\theta_{\widetilde{\ell}_1}^\pm(t)-\pi$. 
  En particulier, $\sum_{t>T_0}\theta_{\widetilde{\ell}_1}^\pm(t)-\pi$ vaut $0$ si chacune de ces valeurs est nulle et $+\infty$ sinon. Puisque $\ell_1([T_0,+\infty[)$ 
n'est pas contenue dans un 
cylindre plat non dégénéré, il existe $t_\pm$ tels que $\theta_{\widetilde{\ell}_1}^\pm(t_\pm)-\pi>0$, donc   
ces deux sommes valent $+\infty$, et donc $\widetilde{\ell}_{1|[0,+\infty[}$ n'est pas quasi-droit. D'après le lemme \ref{Dank10}, les deux géodésiques coïncident
à partir d'un certain temps, à changement d'origine près.\cqfd

\blemm\label{spiralesurcylindre}
Soient $\ell_1:\RR\to\Sigma$ une géodésique locale périodique de $\srfce$ dont l'image est contenue dans un cylindre plat non dégénéré, $\widetilde{\ell}_1$
un relevé de
$\ell_1$ dans $\widetilde{\Sigma}$ et $\widetilde{B}$ la bande plate maximale contenant $\widetilde{\ell}_1(\RR)$, dont on fixe des paramétrages des composantes
de bord. Soit $\widetilde{\ell}_2$ une géodésique de $\revet$ 
telle que $\widetilde{\ell}_2(+\infty)=\widetilde{\ell}_1(+\infty)$. Alors si $\widetilde{\ell}_2(\RR)$ n'est pas contenue dans $\widetilde{B}$, la géodésique 
$\widetilde{\ell}_2$ coïncide avec un 
des bords $\widetilde{b}$ de $\widetilde{B}$ à partir d'un certain temps, quitte à changer les origines, et donc $\ell_{2}=p\circ\widetilde{\ell}_{2}$ coïncide
avec $b=p(\widetilde{b})$ à partir de ce temps. 
\elemm

\dem La démonstration de ce lemme se déduit facilement de la démonstration de \cite[Prop.~3.10]{Dank10}. Quitte à inverser leurs orientations,
les bords de $\widetilde{B}$ ont
les mêmes couples de points à l'infini que $\widetilde{\ell}_1$ et ils n'ont qu'un seul côté quasi-droit. On considère le bord $\widetilde{b}$ de 
$\widetilde{B}$ dont le côté qui n'est pas quasi-droit contient $\widetilde{\ell}_2(\RR)$. Puisque que $\revet$ est $\CAT(0)$, si $\widetilde{b}(\RR)$ et 
$\widetilde{\ell}_2(\RR)$ ne sont pas disjointes, elles sont confondues le long d'un rayon géodésique. Supposons par l'absurde qu'elles sont disjointes.
Il existe une suites de points de
$\widetilde{b}(\RR)$ qui tend vers $\widetilde{b}(+\infty)$, en lesquels l'angle de $\widetilde{b}$ du côté qui n'est pas quasi-droit est strictement supérieur à $\pi$, et donc il existe des géodésiques dont
les 
origines sont ces points, qui ne rencontrent $\widetilde{b}(\RR)$ qu'en un point. De plus, si les origines de ces géodésiques sont assez proches de 
$\widetilde{b}(+\infty)$, d'après le lemme \cite[Prop.~3.3]{Dank10}, leurs deux extrémités dans $\partial_\infty\widetilde{\Sigma}$ sont comprises entre 
$\widetilde{b}(+\infty)$ et $\widetilde{\ell}_2(-\infty)$ (pour l'ordre $o_\infty$ sur $\partial_\infty\widetilde{\Sigma}$), et alors leurs intersections avec 
$\widetilde{\ell}_2(\RR)$ ont deux composantes connexes, ce qui est impossible dans un
espace métrique $\CAT(0)$. Donc $\widetilde{b}(\RR)$ et $\widetilde{\ell}_2(\RR)$ sont confondues le long d'un rayon géodésique.\cqfd

\medskip

De même que dans la partie \ref{partiecomposanteminimale}, on considère  une métrique hyperbolique complète $m$ sur $\Sigma$, 
ainsi qu'une lamination plate $\Lq$ minimale qui n'est pas une paire de feuilles périodiques opposées et $\Lm$ la lamination hyperbolique minimale qui lui correspond.
On note $\Lqr$ et $\Lmr$ leurs images réciproques par $p$.

\blemm\label{spiralesurcomposanteminimalerecurrente}
Supposons que $\Lambda_{[q]}$ soit un ensemble de feuilles de Levitt denses dans un domaine (voir le lemme \ref{composanteminimalerécurrente}).
Alors si $\lambda$ est une géodésique de $(\Sigma,m)$
qui spirale
sur $\Lambda_m$, et si $\ell$ est la géodésique locale plate qui lui correspond (qui est unique à changement d'origine près), il existe une feuille  
$\ell_1$ de $\Lambda_{[q]}$ et telles que $\ell$ et $\ell_1$ coïncident à partir d'un certain temps (quitte à changer d'origines). 
\elemm

\dem D'après le lemme \ref{spirale}, il existe une feuille $\ell_0$ de $\Lq$ dont un relevé à le même point à l'infini qu'un relevé de $\ell$. Puisque les feuilles de $\Lq$
sont des feuilles de Levitt non périodiques, et donc ne sont pas positivement compactes, d'après le lemme \ref{spiralesurcompact}, $\ell$ n'est pas positivement compacte. Puisqu'elle
n'est pas auto-entrelacée, d'après le lemme \ref{directionconstante}, elle ne rencontre pas de singularité et est donc de direction constante à partir d'un certain
temps, soit $T>0$. Comme $\ell_0$ est dense dans le support de $\Lq$, qui est un domaine d'un feuilletage de direction constante (voir le lemme
\ref{composanteminimalerécurrente}), on peut paramétrer $\ell_0$ pour que son origine appartienne à l'intérieur du support de $\Lq$. Alors, d'après le lemme
\ref{composanteminimale}, il existe $t>T$ tel que $\ell(t)$ soit arbitrairement proche de $\ell_0(0)$, et donc appartienne à l'intérieur du support de $\Lq$.
Puisque $\ell$ et $\ell_0$ ne sont pas entrelacées, si $\ell([t,+\infty[)$ n'est pas contenue dans l'image de $\ell_0$, la géodésique $\ell$ n'intersecte pas $\ell_0$, 
et puisque l'image de $\ell_0$ est dense dans le support de $\Lq$, il est de même direction. Donc $\ell_{|[t,+\infty[}$ est un rayon d'une (ou de deux) feuille de
Levitt
de même direction que $\ell_0$, et est donc un rayon d'une feuille de $\Lq$ d'après le lemme \ref{composanteminimalerécurrente}.\cqfd

\section{Lamination plate minimale de support un graphe fini.}\label{graphefini}

Dans cette partie, on cherche à savoir à quelle condition un graphe métrique fini 
peut être plongé isométriquement dans une surface compacte munie d'une structure de demi-translation telle qu'il soit le support d'une lamination plate minimale (qui ne soit pas une paire
de feuilles périodiques opposées). 
Nous utiliserons les définitions et conventions sur les graphes de \cite{Serre83}. On confond un graphe et sa réalisation topologique,
et on le munit d'une distance géodésique telle que chaque demi-arête soit isométrique à un
intervalle compact de $\RR$. On supposera de plus qu'il est connexe et sans sommet terminal. La {\it valence} d'un point est le cardinal de l'ensemble des germes de rayons 
issus de ce point. Un {\it arbre simplicial} est un graphe simplement connexe dont la longueur des arêtes est minorée par une constante strictement positive. 



Soit $G$ un graphe. Une {\it orientation cyclique} de $G$ est la donnée, pour chaque sommet $x$, d'un ordre cyclique sur l'ensemble des germes de
rayons géodésiques issus de $x$.    
  On appelle {\it graphes cycliquement ordonnés} les graphes 
munis d'une orientation cyclique.

\brema\label{ordrecoherent} Soit $(T,or)$ un arbre simplicial cycliquement ordonné, et $\partial_\infty T$ son bord à l'infini.     
             Si $a,b,c$ sont des points deux à deux distincts de $\partial_\infty T$, on note $t=t(a,b,c)$ le centre du tripode défini par $a,b$ et $c$
($t=]a,b[\cap]a,c[\cap]b,c[$), et $r_a,r_b,r_c$ les germes issus de $t$ définis par $[t,a[,[t,b[,[t,c[$. On pose alors $o(a,b,c)=or_{t}(r_a,r_b,r_c)$. La fonction 
$o:(\partial_\infty T)^3\to\{-1,0,1\}$ ainsi définie (avec $o(a,b,c)=0$ si $a,b,c$ ne sont pas deux à deux distincts) est un ordre cyclique sur $\partial_\infty T$, et 
l'application $f$ qui à une orientation cyclique de $T$ associe 
l'ordre cyclique total sur $\ddT$ ainsi défini est injective (voir \cite[Prop.~3.9]{Wolf11}). On dit qu'un ordre cyclique total sur $\ddT$ est {\it cohérent}
s'il appartient à l'image de $f$.  
\erema

Si $(G,or)$ est un graphe fini cycliquement ordonné, et si 
$p:\widetilde{G}\to G$ est un revêtement universel localement isométrique, alors $\widetilde{G}$ est un arbre simplicial et il
existe une unique orientation cyclique $\widetilde{or}$ sur $\widetilde{G}$
telle que $p$ préserve l'orientation cyclique, c'est-à-dire que si $r_1,r_2,r_3$ sont des germes de rayons géodésiques issus d'un
point $\widetilde{x}\in\widetilde{G}$,     
alors 
$or_{x=p(\widetilde{x})}(p(r_1),p(r_2),p(r_3))=\widetilde{or}_{\widetilde{x}}(r_1,r_2,r_3)$. Et cette orientation cyclique définit un ordre cyclique cohérent sur
$\partial_\infty\widetilde{G}$. Un graphe fini cycliquement ordonné est donc aussi un espace métrique enrubanné.

Si $(T,or)$ est un arbre simplicial cycliquement ordonné dont la valence de chaque sommet est finie, en particulier si $T$ est un revêtement d'un graphe fini
cycliquement ordonné, c'est un espace
métrique $\CAT(0)$ propre, et son bord à l'infini est muni d'un ordre (cyclique total). Donc on est dans les conditions de la partie \ref{CAT}, et l'on 
peut définir des laminations géodésiques sur de tels graphes.

 \medskip

Soient $(G,or)$ un graphe cycliquement ordonné et $(\Sigma,[q])$ une surface munie d'une structure de demi-translation tels qu'il existe un plongement isométrique $\phi$ de $G$ dans $\srfce$.
On dit que $\phi$ {\it préserve l'orientation cyclique} si en notant $or^\Sigma$ l'orientation cyclique induite sur le graphe plongé $\phi(G)$ par 
un choix d'orientation de $\Sigma$, alors
 $\phi:(G,or)\to(\phi(G),or^\Sigma)$ préserve l'orientation cyclique.


\smallskip

\noindent
\begin{minipage}{12.5 cm}
Une géodésique locale $g:\RR\to G$ {\it serre à droite}
si pour tout $t\in\RR$ tel que $g(t)$ est un sommet de $G$, le germe de $g([t,t+\varepsilon[)$ succède à celui de $g(]t-\varepsilon,t])$ 
(pour $\varepsilon>0$), pour l'ordre $or_{g(t)}$ (c'est-à-dire, si on note $r_1$ et $r_2$ les germes de $g([t,t+\varepsilon[)$ et de $g(]t-\varepsilon,t])$, 
alors pour tout germe $r_3$ distinct de $r_1$ et $r_2$, on a $or_{g(t)}(r_3,r_1,r_2)=1$). 
 
\end{minipage}    
\begin{minipage}{2.4 cm}
\begin{center}
\begin{picture}(0,0)%
\includegraphics{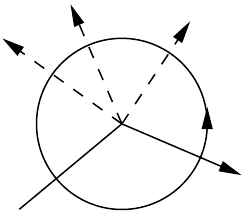}%
\end{picture}%
\setlength{\unitlength}{3771sp}%
\begingroup\makeatletter\ifx\SetFigFont\undefined%
\gdef\SetFigFont#1#2#3#4#5{%
  \reset@font\fontsize{#1}{#2pt}%
  \fontfamily{#3}\fontseries{#4}\fontshape{#5}%
  \selectfont}%
\fi\endgroup%
\begin{picture}(1225,1054)(9246,-235)
\put(9806,-48){\makebox(0,0)[lb]{\smash{{\SetFigFont{8}{9.6}{\familydefault}{\mddefault}{\updefault}{\color[rgb]{0,0,0}$g$}%
}}}}
\end{picture}%

\end{center}
\end{minipage}
\noindent
Chaque arête de $G$ est parcourue dans les 
deux sens par deux géodésiques locales qui serrent à droite, éventuellement égales, qui sont uniques (à changement d'origine près).

\blemm\label{grapheplongé}
Soit $(G,or)$ un graphe fini cycliquement ordonné. Alors $G$ est le support d'une lamination minimale non dénombrable, donc dont aucune feuille n'est positivement ou négativement
périodique,
si et seulement si $(G,or)$ n'est pas isomorphe à un cercle, une paire d'haltères, un huit plat ou un théta plat, par un isomorphisme préservant l'orientation cyclique
(c'est-à-dire \begin{picture}(0,0)%
\includegraphics{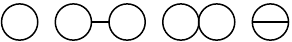}%
\end{picture}%
\setlength{\unitlength}{3771sp}%
\begingroup\makeatletter\ifx\SetFigFont\undefined%
\gdef\SetFigFont#1#2#3#4#5{%
  \reset@font\fontsize{#1}{#2pt}%
  \fontfamily{#3}\fontseries{#4}\fontshape{#5}%
  \selectfont}%
\fi\endgroup%
\begin{picture}(1461,198)(-8,830)
\end{picture}%
, où les orientations cycliques sont données par l'orientation du plan).
\elemm

\dem Décrivons le plongement classique de $G$ dans une surface munie d'une structure de demi-translation. Soit $g:\RR\to G$ une géodésique locale qui serre à droite. Puisque $G$ est fini, il existe au moins deux réels $t_1<t_2$ et $\varepsilon>0$ tels que les germes de 
$g_{|[t_1,t_1+\varepsilon[}$ et de $g_{|[t_2,t_2+\varepsilon[}$ sont égaux, et puisque $g$ est déterminée par la donnée d'un de ses germes de rayon géodésique, 
la géodésique locale $g$ est périodique.
Puisque $G$ est fini, l'ensemble des géodésiques locales qui serrent à droite (à changements d'origine près) est fini et chaque arête est
parcourue dans les deux sens par deux éléments de cet ensemble.


Le long de chaque élément de cet ensemble, on recolle, par isométrie locale, un cylindre plat de hauteur $1$ et de la bonne circonférence, le long 
de l'une de ses composantes de bord et on note $(\Sigma,d)$ l'espace métrique ainsi obtenu et toujours  $G$ le plongé (isométrique) de $G$ dans $(\Sigma,d)$
(voir la fin de la partie pour des exemples).  
Alors par construction, chaque point de $\Sigma$ qui n'est pas un sommet de valence supérieure ou égale à trois de $G$ admet un voisinage isométrique à
un ouvert du demi-plan supérieur fermé de $\RR^2$
(muni de la métrique euclidienne), et chaque sommet de valence supérieure ou égale à trois
de $G$ est une singularité conique d'angle $k\pi$, où $k\geqslant 3$ est sa valence.    
      De plus, les cylindres plats sont recollés le long de 
géodésiques locales de direction
constante. Donc si on les feuillette par des feuilles parallèles à leur bord, on construit ainsi 
un champ de droites parallèle sur $(\Sigma,d)$. Enfin, chacun des cylindres plats est orientable, et puisque les recollements en chaque sommet sont ordonnés selon 
l'ordre cyclique sur les germes de rayons issus du sommet, la surface $\Sigma$ est orientable. On a vu dans la partie \ref{diffquadholo} que $\Sigma$ est 
munie d'une structure de demi-translation, soit $[q]$. 
On peut choisir l'orientation de $\Sigma$ pour que le plongement de $G$ dans $\Sigma$ préserve l'orientation.


La surface compacte $\Sigma$ ne peut être ni un disque, ni une sphère, ni un tore. Supposons que $\Sigma$ n'est ni un cylindre ni un pantalon. Alors
$\chi(\Sigma)<0$ et il existe une métrique hyperbolique $m$ à bord totalement géodésique sur $\Sigma$. D'après le lemme \ref{maximaleminimale},
il existe au moins une 
lamination minimale et quasi-maximale sur $\srfcem$, qui est non dénombrable. Soit  $\Lq$ l'unique lamination plate qui lui correspond, qui est aussi non dénombrable
(voir la partie \ref{hyperboliqueplate} pour la
 définition de cette correspondance). On note $\Lqr$ et $\Lmr$ les images réciproques de $\Lq$ et $\Lm$ dans $\widetilde{\Sigma}$.  Puisque $\Lm$ est minimale 
 et n'est pas réduite à une feuille fermée, d'après le lemme \ref{minimalepaspériodique}, aucune feuille de $\Lq$ 
 n'est positivement ou négativement périodique.  

Supposons qu'il existe un point d'une feuille $\ell$ de $\Lq$ qui n'appartienne pas à $G$. Si $\ell$ n'est pas parallèle à une composante du bord, il existe un rayon de
$\ell$ qui rencontre transversalement une composante du bord, ce qui est impossible car alors $\ell$ ne serait pas définie sur $\RR$. Donc $\ell$ est parallèle à 
une composante du bord, et donc périodique, ce qui est impossible car $\Lq$ est minimale et non dénombrable. Donc le support de $\Lq$ est contenu dans $G$.

 \begin{center}
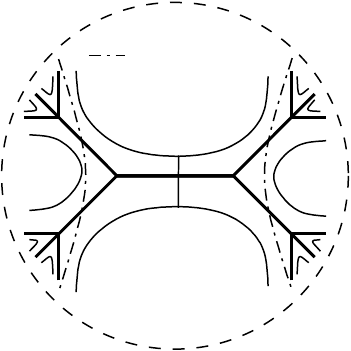
\end{center}
 Supposons qu'il existe une arête $a$ de $G$ qui n'est parcourue par aucune feuille de $\Lq$. Alors si $\widetilde{a}:[-r,r]\to\widetilde{\Sigma}$ est un relevé paramétré
 de $a$ (avec $r>0$),
 en notant $A$ et $B$ les ensembles de points de $\partial_\infty\widetilde{\Sigma}$ qui sont respectivement les extrémités des rayons qui ont le même germe que
 $\widetilde{a}_{|[0,r[}$ ou $\widetilde{a}_{|]-r,0]}$, aucune feuille de $\Lqr$ n'a ses points à l'infini appartenant respectivement à $A$ et $B$, et de même pour $\Lmr$.
 Donc le plus court segment géodésique de $\revetm$ qui joint les composantes de bord $\widetilde{b}_1$ et $\widetilde{b}_2$ (voir la  figure), n'est
 intersecté par aucune feuille de $\Lmr$, 
et il existe une composante du complémentaire du support de $\Lmr$ dont le complété a au moins $\widetilde{b}_1$ et $\widetilde{b}_2$ comme composantes de 
$\partial\widetilde{\Sigma}$ dans son bord. Or, puisque $\Lm$ est quasi-maximale, les complétés des composantes
connexes du complémentaire de 
$\widetilde{\Lambda}_{\widetilde{m}}$ dans $\widetilde{\Sigma}$ sont des triangles géodésiques idéaux ou des polygones géodésiques idéaux dont un seul côté n'est pas une feuille mais
une composante du bord de $\widetilde{\Sigma}$. Donc chaque arête de $G$ est parcourue par au moins
une feuille de $\Lq$. D'après le lemme 
\ref{composanteminimalecompacte}, chaque feuille parcourt toutes les arêtes.

\medskip

Si $\Sigma$ est un cylindre plat, les feuilles d'une lamination de $\srfce$ sont périodiques. De même, si $\Sigma$ était un pantalon et s'il existait
une
lamination non dénombrable sur $\srfce$, la lamination hyperbolique lui correspondant serait non dénombrable, ce qui est impossible sur un pantalon hyperbolique. Il reste donc à 
déterminer les graphes finis enrubannés dont la surface munie d'une structure de demi-translation associée, construite ci-dessus, soit un pantalon. Soient $(G,or)$ un tel graphe
et $\srfce$ le pantalon plat associé. Puisque $\Sigma$ se rétracte par déformation forte sur $G$,
le groupe 
fondamental de $G$ en un point quelconque est isomorphe au groupe fondamental de $\Sigma$ au même point, et est donc un groupe libre à deux générateurs.
Le graphe $G$ est donc homéomorphe à un graphe de type huit, haltères, ou  théta. 
Si $G$ est homéomorphe à un graphe de type haltères, quitte à prendre le symétrique d'une haltère par rapport à l'axe,
on montre qu'il n'existe qu'une seule orientation cyclique sur $G$, sinon il existe deux
orientations cycliques sur $G$ (à isomorphisme préservant l'orientation cyclique près). On montre alors que $\Sigma$ n'est un pantalon que si $G$ est le graphe
des haltères, du théta plat ou du huit plat. \cqfd

\begin{center}
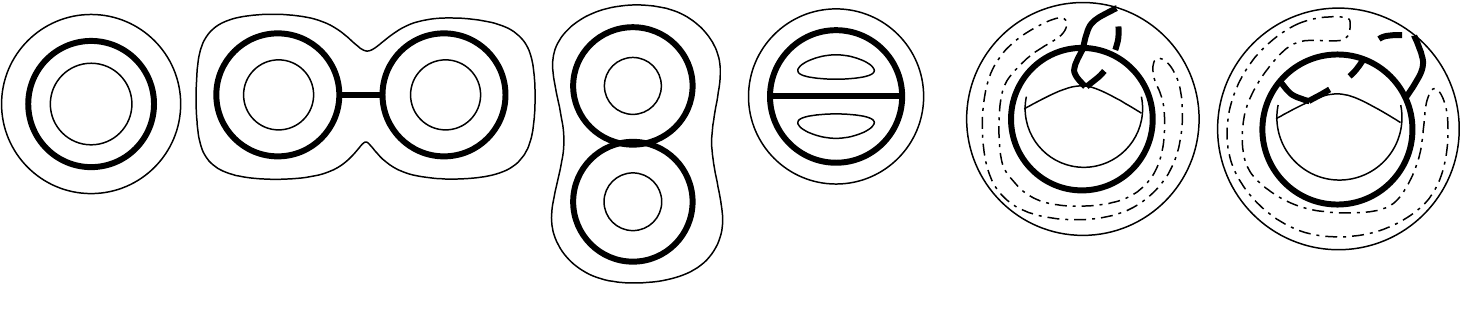
\end{center}

\section{Théorème de structure des laminations géodésiques plates.}\label{laminationplate}

Dans cette partie, on donne un théorème de structure des laminations géodésiques sur une surface compacte munie d'une structure de demi-translation.
On note de nouveau $(\Sigma,[q])$ une surface connexe, orientable,  à bord (éventuellement vide), munie d'une 
structure de demi-translation et on 
suppose que $\Sigma$ est compacte et $\chi(\Sigma)<0$.

\medskip

Soit $\Lq$ une lamination plate. On appelle {\it composante cylindrique} de $\Lq$ un ensemble maximal de feuilles de $\Lq$ contenues dans un cylindre plat 
non dégénéré (elles sont alors périodiques)
et on dit qu'elle est {\it pleine} si son support est égal au cylindre. Une {\it composante minimale} de $\Lq$ est une sous-lamination plate minimale de $\Lq$. 
Une composante minimale est dite de {\it type graphe fini} si son support est un graphe fini (nécessairement sans sommet terminal),
et de type {\it récurrent} si c'est la réunion des feuilles de Levitt locales contenues dans son support et de même direction (voir la définition précédant
le lemme \ref{composanteminimalerécurrente}). 
    Si $\ell$ est une géodésique plate, on dit qu'un de ses bouts 
{\it aboutit} dans une composante cylindrique s'il existe un rayon géodésique correspondant à ce bout qui coïncide avec un bord du cylindre plat maximal
correspondant, à partir d'un certain temps (elle peut  être disjointe des feuilles de la composante cylindrique), et qu'elle {\it aboutit}
dans une composante minimale s'il existe un rayon géodésique correspondant à ce  bout qui est un rayon d'une feuille de cette composante minimale.

\btheo\label{structurelaminationplate}
Soit $\Lq$ une lamination plate sur une surface compacte munie d'une structure de demi-translation $\srfce$. Alors $\Lq$
est la réunion d'un nombre fini de composantes cylindriques, d'un nombre fini de
composantes minimales (de type récurrent ou graphe fini ou  paire de feuilles périodiques opposées d'images non contenues dans un cylindre plat non dégénéré)
et d'un nombre fini de feuilles isolées (pour la topologie des géodésiques) dont chacun des
bouts aboutit dans une composante minimale ou dans une composante cylindrique. 

\etheo


\rem  Les supports des composantes minimales de type récurrent sont des domaines de
feuilletages verticaux de différentielles quadratiques holomorphes (voir le lemme \ref{structuredestrajectoires}).
Si les supports de deux composantes minimales de type récurrent s'intersectent, puisque tous les points de leur support appartiennent à des feuilles de Levitt deux à
deux non entrelacées,
elles 
ne peuvent s'intersecter qu'en une réunion de points isolés et de liaisons de singularités qui n'est pas connexe a priori, le long de leur bord.
De même, une composante
minimale de type récurrent ne peut intersecter une composante cylindrique que le long du bord du cylindre plat maximal qui la contient, 
et deux cylindres plats maximaux ne peuvent s'intersecter que le long de leur bord. On ne peut rien dire, a priori de l'intersection des autres 
types de composantes.

\medskip

\dem On fixe une métrique hyperbolique $m$ (complète à bord totalement géodésique) sur $\Sigma$.
On renvoie à  la partie \ref{hyperboliqueplate} pour les notations et les définitions de lamination plate pleine et de lamination hyperbolique correspondant 
à une lamination plate.  
 Soit $\Lm$ la lamination géodésique hyperbolique qui correspond à $\Lq$. Puisque $\Sigma$ est compacte, d'après par exemple \cite[Th.~1.4.2.8~p.~83]{CaEpMa05}, la lamination hyperbolique $\Lm$ est 
la réunion disjointe d'un nombre fini de composantes minimales et d'un nombre fini de géodésiques isolées dont chacun des bouts 
spirale 
sur une composante minimale.
D'après le corollaire \ref{corogeodasymptotes}, à chaque composante minimale de $\Lm$ qui n'est pas une feuille fermée correspond une unique composante minimale de $\Lq$,
qui est de type graphe fini si elle contient une géodésique d'image compacte, et de type récurrent sinon (d'après les lemmes  \ref{composanteminimalecompacte} et 
\ref{composanteminimalerécurrente}). De plus, si une géodésique de $\Lm$ est périodique et si elle ne correspond pas à une unique
feuille 
périodique de $\Lq$, elle correspond à un ensemble maximal
d'au moins deux feuilles de $\Lq$, contenues dans un cylindre plat d'après le lemme \ref{geodasymptotes}, qui est fermé pour la topologie des géodésiques, d'après le lemme \ref{Xfermé}. 
Cet ensemble est donc une composante
cylindrique.

Enfin, si le bout d'une feuille de $\Lm$ spirale sur une composante minimale de $\Lm$, le bout de la géodésique plate qui lui correspond
aboutit sur la composante cylindrique ou la composante minimale de $\Lq$ correspondante, d'après les lemmes \ref{spiralesurcompact}, 
\ref{spiralesurcylindre} et \ref{spiralesurcomposanteminimalerecurrente}.\cqfd

\bibliographystyle{alphanum}
\bibliography{biblio}{}

Département de mathématique, UMR 8628 CNRS, Université Paris-Sud, Bât. 430, F-91405 Orsay Cedex, FRANCE. 
Bureau : 16.

{\it thomas.morzadec@math.u-psud.fr}
\end{document}